\documentclass[a4paper,11pt]{article} 

\title{A Survey on the  Eigenvalues  Local Behavior of \\Large Complex Correlated Wishart Matrices}

\author{
\ Walid Hachem
\footnote{
CNRS LTCI; T\'{e}l\'{e}com ParisTech, 46 rue Barrault, 75046 Paris Cedex 13, France. Email: walid.hachem@telecom-paristech.fr}\;\;,\;
\ Adrien Hardy  
\footnote{Department of Mathematics, KTH Royal Institute of Technology, Lindstedtsv\"agen 25, 10044 Stockholm, Sweden. Email: ahardy@kth.se}\;\;,\; 
\ Jamal Najim 
\footnote{CNRS LIGM; Universit\'e Paris-Est, 
Cit\'e Descartes, 
5 Boulevard Descartes, 
Champs sur Marne, 
77 454 Marne la Vall\'ee Cedex 2,
France. Email: najim@univ-mlv.fr}}

\usepackage[english]{babel}
\usepackage[utf8]{inputenc}

\usepackage{a4wide}
\usepackage{fancyhdr}
\usepackage{amssymb,amsmath,amscd,amsfonts,amsthm,bbm}
\usepackage{tikz,graphicx,subfigure,overpic}
\usepackage{color} 
\usepackage[colorlinks=true, linkcolor=blue, citecolor=blue]{hyperref}
\usepackage[fixlanguage]{babelbib}
\usepackage{color}
\usepackage{hyperref}
\usepackage{comment}

\numberwithin{equation}{section}

\def\Tr{\mathop{\mathrm{Tr}}\nolimits}
\def\re{\mathop{\mathrm{Re}}\nolimits}

\newtheorem{theorem}{Theorem}[]
\newtheorem{lemma}{Lemma}[section]
\newtheorem*{theorem*}{Theorem}

\newtheorem{proposition}[lemma]{Proposition}
\theoremstyle{definition} 
\newtheorem{definition}[lemma]{Definition}
\newtheorem{assumption}{Assumption}[]
\newtheorem{Remark}[lemma]{Remark}
\newenvironment{remark}{\begin{Remark}\rm}{\end{Remark}}

\newtheorem{Example}[lemma]{Example}

\newcommand{\eq}{\begin{equation}}
\newcommand{\qe}{\end{equation}}

\newcommand{\R}{\mathbb{R}}
\newcommand{\C}{\mathbb{C}}

\newcommand{\p}{\mathbb{P}}
\newcommand{\K}{ {\rm K}}

\newcommand{\E}{{\rm E}}

\DeclareMathOperator{\supp}{Supp}

\newcommand{\Ai}{{\rm Airy}}
\newcommand{\Be}{{\rm Bessel}} 
\newcommand{\Pe}{{\rm Pearcey}}

\newcommand{\bt}{\textbf}
\newcommand{\bs}{\boldsymbol}
\newcommand{\bv}{\mathbf}
 
\renewcommand{\leq}{\leqslant}
\renewcommand{\geq}{\geqslant}

\renewcommand{\epsilon}{ \varepsilon}
\renewcommand{\phi}{ \varphi}
\renewcommand{\d}{ {\rm d}}
\renewcommand{\frak}{\mathfrak} 
\renewcommand{\emptyset}{\varnothing}

\begin{document}
\maketitle

\begin{abstract} 

The aim of this note is to provide a pedagogical survey of the recent works
\cite{HHN,HHN2} concerning the local behavior of the eigenvalues of large
complex correlated Wishart matrices at the edges and cusp points of the spectrum: Under quite general conditions, 
the eigenvalues fluctuations at a soft edge of the 
limiting spectrum, at the hard edge when it is present, or at 
a cusp point, are  respectively described by mean of the Airy kernel, the Bessel kernel, or the Pearcey  
kernel. 
Moreover, the eigenvalues fluctuations at several soft edges are asymptotically independent. 
In particular, the asymptotic fluctuations of the matrix condition number 
can be described. Finally, the next order term of the hard edge 
asymptotics is provided.
\end{abstract}

\setcounter{tocdepth}{2}


\section{The matrix model and assumptions} 

Consider the $N\times N$ random matrix defined as
\eq
\label{main matrix model}
{\bv M}_N = \frac 1N {\bf X}_N {\bf \Sigma}_N{\bf X}_N^* 
\qe
where  ${\bf X}_N$ is an $N\times n$ matrix with independent and identically
distributed (i.i.d.) entries with zero mean and unit variance, and 
${\bf \Sigma}_N$ is a $n\times n$ deterministic positive definite Hermitian 
matrix.
The random matrix $\bv M_N$ has $N$ non-negative eigenvalues, but which may be
of different nature. Indeed, the smallest $N-\min(n,N)$ eigenvalues are
deterministic and all equal to zero, whereas the other $\min(n,N)$ eigenvalues
are random. The problem is then to describe the asymptotic behavior of the
random eigenvalues of $\bv M_N$, as both dimensions of ${\bf X}_N$ grow to
infinity at the same rate. 
Let us mention that the $n\times n$ random covariance matrix  
\[
\widetilde{\bv M}_N = \frac 1N 
{\bf \Sigma}_N^{1/2} {\bf X}_N^* {\bf X}_N {\bf \Sigma}_N^{1/2}\ ,
\]
which is also under consideration, has exactly
the same random eigenvalues as $\bv M_N$, and hence  results on the random
eigenvalues can be carried out from one model to the other immediately. 

The global behavior of the spectral distribution of $\widetilde{\bv M}_N$ in 
the large dimensional regime is known since the work of Mar\v cenko 
and Pastur~\cite{MP}, where it is shown that this spectral distribution
converges to a deterministic probability measure $\mu$ that can be identified. 
In this paper, we will be interested in the local behavior of the eigenvalues 
of $\widetilde{\bv M}_N$ near the edge points and near the so-called cusp 
points of the support of $\mu$. The former will be called the extremal 
eigenvalues of $\widetilde{\bv M}_N$. 


The random matrices ${\bv M}_N$ and $\widetilde{\bv M}_N$ are ubiquitous 
in multivariate statistics~\cite{book-bai-chen-liang-2009},  mathematical finance
\cite{bouchaud-noise-dressing-99,guhr-credit-risk-2014},  electrical engineering and signal processing \cite{book-couillet-debbah}, etc. Indeed, in multivariate statistics, the performance 
study of the Principal Component Analysis algorithms \cite{Joh01} requires the 
knowledge of the fluctuations of the extremal eigenvalues of $\widetilde{\bv M}_N$. In mathematical finance, $\widetilde{\bv M}_N$ represents the empirical covariance matrix obtained from a sequence of asset returns.  In signal processing,  $\widetilde{\bv M}_N$ often stands for the empirical
covariance matrix of a spatially correlated signal received by an array of antennas, and source detection \cite{BDMN11,KN09} 
or subspace separation \cite{MDCM} algorithms also rely on 
the statistical study of these extremal eigenvalues.

In this article, except when stated otherwise, we restrict ourselves to the case of complex Wishart matrices. Namely, we make the following assumption. 

\begin{assumption} 
\label{ass:gauss} 
The entries of $\bv X_N$ are i.i.d. standard complex Gaussian random variables.
\end{assumption}

Concerning the asymptotic regime of interest, we  consider here the large random matrix regime, where the number of rows and columns of $\bv M_N$ both grow to infinity at the same pace. More precisely,  we assume $n=n(N)$ and $n,N\to \infty$  in such a way that 
\eq
\label{evdistrMN}
\lim_{N\rightarrow\infty} \frac{n}{N}=\gamma\in (0,\infty)\ .
\qe 
This regime will be simply referred to as  $N\to \infty$ in the sequel.  

Turning to $\bv \Sigma_N$, let us denote by 
$0<\lambda_1\leq \cdots\leq \lambda_n$ the eigenvalues of this matrix and let
\eq
\label{nuN}
\nu_N=\frac 1n \sum_{j=1}^n \delta_{\lambda_j} 
\qe 
be its spectral measure. Then we make the following assumption. 

\begin{assumption}\  
\label{ass:nu}
\begin{enumerate}
\item 
The measure $\nu_N$ weakly converges towards a limiting probability measure
$\nu$ as $N\rightarrow\infty$, namely
\[
\frac{1}{n}\sum_{j=1}^n f(\lambda_j)\xrightarrow[N\to\infty]{} \int f(x)\nu(\d x)
\] 
for every bounded and continuous function $f$.
\item
For $N$ large enough, the eigenvalues of $\bv \Sigma_N$ stay in a compact 
subset of $(0,+\infty)$ independent of $N$, i.e.
\eq
0\ <\ \liminf_{N\rightarrow\infty}\lambda_1,\qquad 
\sup_{N}\lambda_n\ <\ +\infty. 
\qe
\end{enumerate}
\end{assumption}

Under these assumptions, a comprehensive description of the large $N$ behavior
of the eigenvalues of $\bv M_N$ can be made. To start with, we recall in
Section~\ref{section global} some classical results describing the global
asymptotic behavior of these eigenvalues, as a necessary step for studying
their local behavior. We review the results of Mar\v cenko-Pastur~\cite{MP} and
those of 
Silverstein-Choi~\cite{SC}, which show among other things that the spectral
measure of $\bv M_N$ converges to a limit probability measure $\mu$, that $\mu$
has a density away from zero, that the support of $\mu$ can be delineated, and
that the behavior of the density of $\mu$ near the positive endpoints (soft
edges) of this support can be characterized.  We moreover complete the picture
by describing the behavior of the limiting density near the origin when it is
positive there (hard edge), and also when it vanishes in the interior of the
support (cusp point). The latter results are extracted from~\cite{HHN2}. 

Next, in Section \ref{section local} we turn to the eigenvalues local
behavior. More precisely, we investigate the behavior of the random eigenvalues
after zooming around several points of interest in the support, namely the soft
edges, the hard edge when existing, and the cusp points. In a word, it is shown
in the works~\cite{HHN,HHN2} that the Airy
kernel, the Bessel kernel, and the Pearcey kernel describe the local statistics
around the soft edges, the hard edge, and the cusp points respectively,
provided that a regularity condition holds true.  In particular, the extremal
eigenvalues fluctuate according to Tracy-Widom laws.

In Section \ref{sec:proofs}, we provide sketches of proofs. We first recall an important expression for the kernel $\K_N$ associated 
to the (random) eigenvalues of ${\bf M}_N$ and then outline how one can prove asymptotic convergence towards Airy, Pearcey or Bessel kernels by zooming around the points of interest: either a soft edge, a cusp point or the hard edge.

In Section \ref{OP}, we provide a list of open questions, directly related to the results of the paper.

\paragraph*{Acknowledgements.} 
WH is pleased to thank the organizers of the {\em Journ\'ees MAS 2014} where the project of this note was initiated.
During this work, AH was supported by the grant KAW 2010.0063 from the
Knut and Alice Wallenberg Foundation. 
The work of WH and JN was partially supported by the program 
``mod\`eles num\'eriques'' of the French Agence Nationale de la Recherche 
under the grant ANR-12-MONU-0003 (project DIONISOS).

\section{Global behavior}
\label{section global}

Since the seminal work of Mar\v cenko and Pastur~\cite{MP}, it is known that
under Assumptions~\ref{ass:gauss} and~\ref{ass:nu} the spectral measure of 
${\bv M}_N$ almost surely (a.s.) converges weakly towards a 
limiting probability measure $\mu$ with a compact support. Namely we have 
\eq
\label{mu}
\frac{1}{ N}\Tr f({\bv M}_N)\xrightarrow[N\to\infty]{a.s.} 
\int f(x)\mu(\d x)
\qe
for every bounded and continuous function $f$. As a probability measure, $\mu$
can be characterized by its Cauchy transform: this is the holomorphic function 
defined by 
\[
m(z)=\int \frac{1}{z-\lambda}\,\mu(\d \lambda),\qquad 
z\in \C_+ =\big\{z\in\C : \; {\rm Im}(z)>0\big\},
\]
and which takes its values in  $\C_-=\{z\in\C : \; {\rm Im}(z)<0\}$. More precisely, for any open interval $I\subset\R$ with neither endpoints on an atom of $\mu$, we have the inversion formula
\[
 \mu(I)=-\frac{1}{\pi} \lim_{\epsilon\to 0}\int_I {\rm Im} \big(m(x+i\epsilon)\big)\d x .
\]
For every $z\in\C_+$, the Cauchy transform $m(z)$ of $\mu$ happens to be the unique solution 
$m \in \C_-$ of the fixed-point equation
 \begin{equation}
 \label{cauchy eq} 
m = \left( z - \gamma \int \frac{\lambda}{1 - m \lambda} \nu(\d\lambda) 
\right)^{-1} \ , 
 \end{equation}
where $\gamma$ and $\nu$ were introduced in \eqref{evdistrMN} and Assumption~\ref{ass:nu}.

Moreover, using the free probability terminology \cite{AGZ,HP}, the limiting distribution $\mu$ is also known to be the free multiplicative convolution of the Mar\v cenko-Pastur law \eqref{mp}--\eqref{dens-mp}  with $\nu$, and equation \eqref{cauchy eq} is a consequence from the subordination property of the multiplicative free convolution \cite{Cap}.

For example, in the case where $\nu = \delta_1$, which happens e.g. when $\bs\Sigma_N = I_n$, this equation
has an explicit solution and the  measure $\mu$ can be recovered explicitly. 
\eq
\label{mp}
\mu(\d x)=\left(1-\gamma\right)^+ \, \delta_0 + \rho(x)\d x,
\qe
where $x^+=\max(x,0)$ and the density $\rho$ has the expression
\begin{equation}
\label{dens-mp} 
\rho(x)=\frac{1}{2\pi x}\sqrt{(\frak b-x)\big(x-\frak a)}\;\bs 1_{[\frak a,\frak b]}(x),\qquad \frak a=(1-\sqrt\gamma)^2,\qquad \frak b=(1+\sqrt\gamma)^2.
\end{equation} 
This is the celebrated Mar\v cenko-Pastur law. 

\begin{figure}[h]
\centering
\includegraphics[width=0.7\linewidth]{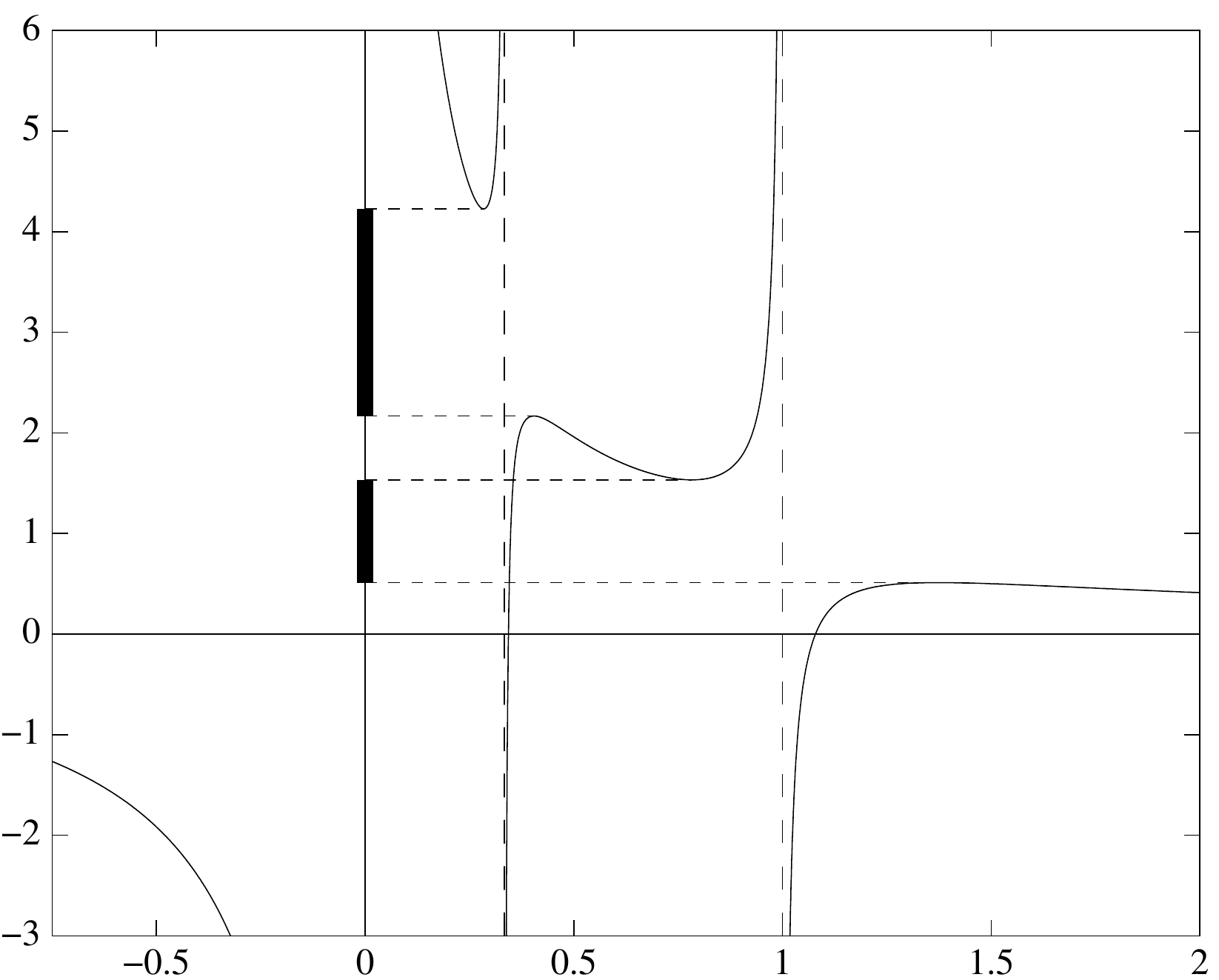}
\caption{Plot of $g:D\to \mathbb{R}$ for $\gamma=0.1$ and 
$\nu = 0.7\delta_1 + 0.3\delta_3$. In this case, 
$D~=~(-\infty,0)~\cup~(0,\frac 13)~\cup~ (\frac 13,1)~\cup~(1,\infty)$. 
The two thick segments on the vertical axis represent $\supp(\rho)$.}
\label{fig:g2bulks}
\end{figure}

\begin{figure}[h]
\centering
\includegraphics[width=0.7\linewidth]{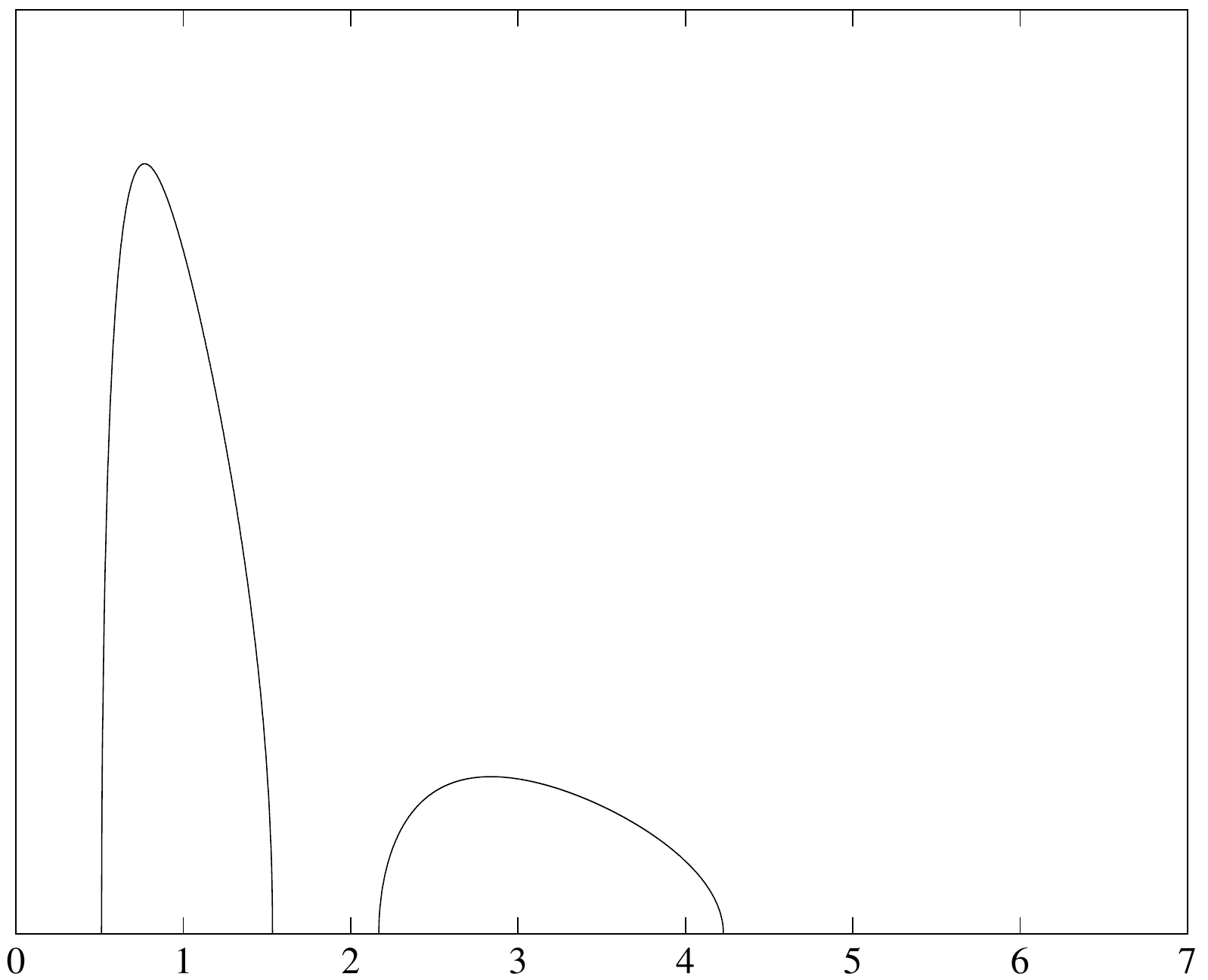}
\caption{Plot of the density $\rho$ in the framework of 
Figure~\ref{fig:g2bulks}.} 
\label{fig:f2bulks}
\end{figure}

When $\nu$ has a more complicated form, it is in general impossible to obtain an
explicit expression for $\mu$, except in a few particular cases.  Nonetheless,
it is possible to make a detailed analysis of the properties of this measure,
and this analysis was done by Silverstein and Choi in~\cite{SC}. These 
authors started by
showing that $\lim_{z\in\C_+ \to x} m(z) \equiv m(x)$ exists for every 
$x \in \R^* = \R - \{ 0 \}$. Consequently, the function $m(z)$ can be 
continuously extended to $\C_+ \cup \R^*$, and furthermore, $\mu$ has a 
density on $\R^*$ defined as $\rho(x) = - \pi^{-1} {\rm Im} \, (m(x))$. We still have the representation  
\eq
\label{densite} 
\mu(\d x)=\left(1-\gamma\right)^+ \, \delta_0 + \rho(x)\d x
\qe
 with this new $\rho$,  making $\rho(x) \d x$  the limiting distribution of the random eigenvalues 
of $\bv M_N$. As is common in random matrix theory, we shall refer to 
the support of $\rho(x)\d x$ as the \textbf{bulk}; we will denote (with a slight abuse of notation) its support by $\supp(\rho)$. 
Silverstein and Choi
also showed that $\rho$ is real analytic wherever it is positive, and they
moreover characterized the compact support $\supp(\mu)$ following the ideas
of~\cite{MP}. More specifically, one can see that the function $m(z)$
has an explicit inverse (for the composition law) on $m(\C_+)$ defined by
\begin{equation}
\label{g(m)} 
g(m) = \frac{1}{m} + 
    \gamma \int \frac{\lambda}{1 - m \lambda} \nu(\d\lambda) , 
\end{equation}
and that this inverse extends to $\C_- \cup D$ and is real analytic on $D$, 
where $D$ is the open subset of the real line
\begin{equation}
\label{D}
D=\bigl\{x\in\R : \ x\neq 0, \ x^{-1}\notin \supp(\nu)\bigr\} . 
\end{equation} 
It was proved in \cite{SC} that 
\[
\R - \supp(\rho) = \big\{ g(m) : \;  m\in D, \ g'(m) < 0 \big\} .
\]
An illustration of these results is provided by Figures~\ref{fig:g2bulks} 
and~\ref{fig:f2bulks}. 

Of interest in this paper are the left edges, the right edges and the cusp 
points of $\supp(\rho)$. 

\noindent
A \textbf{left edge} is a real number $\frak a$ satisfying for every $\delta>0$ small enough
\[
\int_{\frak a-\delta}^\frak{a}\rho( x)\d x=0,\qquad \int_{\frak a}^\frak{a+\delta}\rho( x)\d x>0\ .
\]
A  \textbf{right edge} is a real number $\frak a$ satisfying for every $\delta>0$ small enough
\[
\int_{\frak a-\delta}^\frak{a}\rho(x)\d x>0,\qquad \int_{\frak a}^\frak{a+\delta}\rho( x)\d x=0\ .
\]
A \textbf{cusp point} is a real number $\frak a$ such that $\rho(\frak a)=0$ and, for every $\delta>0$ small enough, 
\[
 \int_{\frak a-\delta}^\frak{a}\rho(x)\d x>0\quad \textrm{and}\quad \int_{\frak a}^\frak{a+\delta}\rho(x)\d x>0\ .
\]

Of course all edges and cusp points are positive numbers, except perhaps the
leftmost edge. When the leftmost edge is the origin, it is common in random
matrix theory to refer to it as the \textbf{hard edge}. In contrast, any
positive edge is also called a \textbf{soft edge}.

The results of~\cite{SC} summarized above show that the study of the map $g$ on
the closure $\overline D$ of $D$ provides a complete description for the edges
and the cusp points.  First, a right edge is either a local minimum of $g$
reached in $D$, or belongs to $g(\partial D)$, which means there is $\frak c
\in\partial D=\overline D\setminus D$ such that 
$\lim_{x\to\frak c,\, x\in D}g(x)$ 
exists, is finite, and equals to that edge. In the former case, $\rho(x)$ 
behaves like a square root near the edge. 
 
\begin{proposition}
\label{prop global R}
If   $\frak a$ is a right edge, then either there is a unique $\frak c\in D$ such that 
\eq
\label{g right}
g(\frak c)=\frak a, \qquad g'(\frak c)=0,\qquad g''(\frak c)>0,
\qe
or  $\frak a\in g(\partial D)$. In  the former case,  we have 
\eq
\label{sqrt right}
\rho(x)=\frac{1}{\pi}\left(\,\frac{2}{g''(\frak c)}\,\right)^{1/2}\big(\frak a-x\big)^{1/2}\;(1+o(1))\ , \qquad x\to\frak a_- \ .
\qe
Conversely , if $\frak c\in D$ satisfies \eqref{g right}, then $\frak a$ is a right edge and \eqref{sqrt right}  holds true.
\end{proposition}

The case where an edge lies in $g(\partial D)$ turns out to be quite delicate.
In the forthcoming description of the eigenvalues local behavior  near the
edges, we shall restrict ourselves to the edges arising as local minima of $g$,
see also Section \ref{OP} for further discussion. Notice also that if $\nu$ is
a discrete measure, as exemplified by Figures~\ref{fig:g2bulks} 
and~\ref{fig:f2bulks}, then $g$ is infinite on $\partial D$ and in particular a right edge cannot belong to $g(\partial D)$:  
the right edges are in this case in a one-to-one correspondence with the local
minima of $g$ on $D$.  

The situation is similar for the soft left edges, except that they correspond 
to local maxima.

\begin{proposition}
\label{prop global L}
If   $\frak a>0$ is a left edge, then either there is a unique $\frak c\in D$ such that 
\eq
\label{g left}
g(\frak c)=\frak a, \qquad g'(\frak c)=0,\qquad g''(\frak c)<0,
\qe
or  $\frak a\in g(\partial D)$. In  the former case,  we have 
\eq
\label{sqrt left}
\rho(x)=\frac{1}{\pi}\left(\frac{2}{-g''(\frak c)}\right)^{1/2}\big(x-\frak a\big)^{1/2}\;(1+o(1))\ , \qquad x\to\frak a_+\ .
\qe
Conversely, if $\frak c\in D$ satisfies \eqref{g left}, then $\frak a$ is a right edge and \eqref{sqrt left}  holds true.
\end{proposition}
Propositions \ref{prop global R} and \ref{prop global L} have been established in \cite{SC}. We state below their counterparts for the hard edge and a cusp point.

The hard edge setting turns out to be similar to the  soft left edge one, except that $\frak c$ is now located at infinity, and  $\rho(x)$ behaves like an inverse square root near the hard edge. More precisely, observe that the map $g$ is holomorphic at $\infty$ and $g(\infty)=0$, in the sense that the map $z\mapsto g(1/z)$ is holomorphic at zero and vanishes at $z=0$. We also denote by $g'(\infty)$ and $g''(\infty)$ the first and second derivatives of the latter map evaluated at $z=0$.

\begin{proposition}
\label{prop global 0}  
The bulk presents a hard edge if and only if 
\eq
\label{g hard}
g(\infty)=0,\qquad g'(\infty)=0,\qquad g''(\infty)<0,
\qe
or equivalently if $\gamma=1$. In this case, we have 
\eq
\label{sqrt hard}
\rho(x)=\frac{1}{\pi}\left(\frac{2}{-g''(\infty)}\right)^{-1/2}  x^{-1/2}\;(1+o(1))\ ,\qquad x\to 0_+\ .
\qe
\end{proposition}

More precisely, we have the explicit formulas $g'(\infty)=1-\gamma$ and
$g''(\infty)=-2\gamma \int \lambda^{-1} \nu(\d \lambda)$. In particular the
statement $g''(\infty)<0$ is always true, and so is $g(\infty)=0$ as explained
above; we included them in \eqref{g hard} to stress the analogy 
with~\eqref{g left}.

A simple illustration of Propositions~\ref{prop global R} 
to~\ref{prop global 0} is provided by the Mar\v cenko-Pastur law. 
From~\eqref{dens-mp}, one immediately sees that 
$\rho(x)\sim  (\frak b- x)^{1/2}$ as $x\to \frak b_-$, and that a
similar square root behavior near $\frak a$  holds if and only if $\frak a> 0$,
that is $\gamma\neq 1$. If $\gamma=1$, i.e., $\frak a=0$, then 
$\rho(x)\sim x^{-1/2}$ as $x\to 0_+$ instead.

We  now turn to the cusp points. Those who will be of interest here correspond
to inflexion points of $g$ where this function is non decreasing. Moreover, a cubic root behavior for the density $\rho(x)$ is observed near such a cusp point, hence justifying the terminology (we recall cusp usually refers to the curve defined by $y^2=x^3$).

\begin{figure}[h]
\centering
\includegraphics[width=0.7\linewidth]{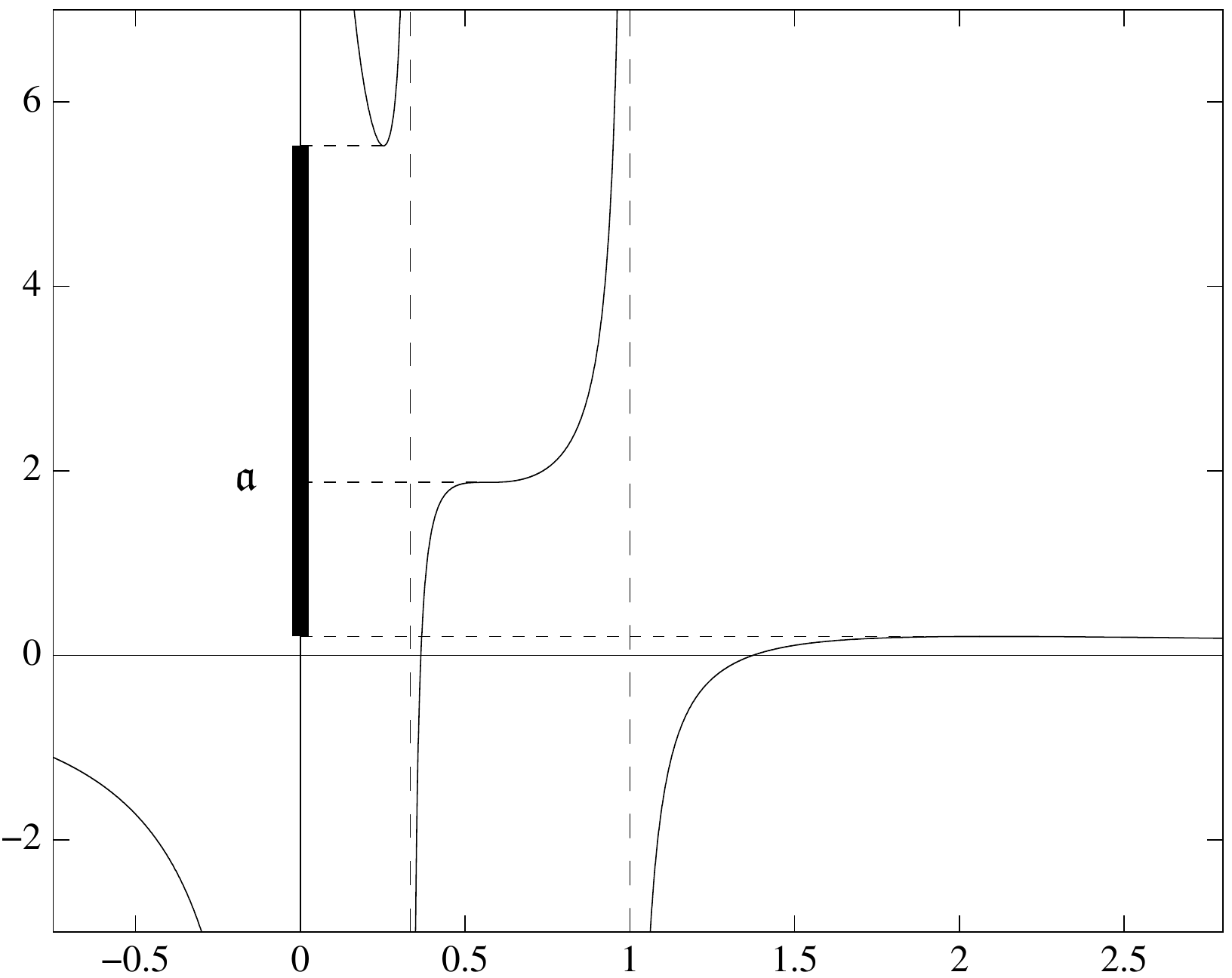}
\caption{Plot of $g$ for $\gamma\simeq 0.336$ and 
$\nu = 0.7\delta_1 + 0.3\delta_3$. The thick segment on the vertical axis 
represents $\supp(\mu)$. The point $\frak a$ is a cusp point. 
}
\label{fig:gcusp}
\end{figure}

\begin{figure}[h]
\centering
\includegraphics[width=0.7\linewidth]{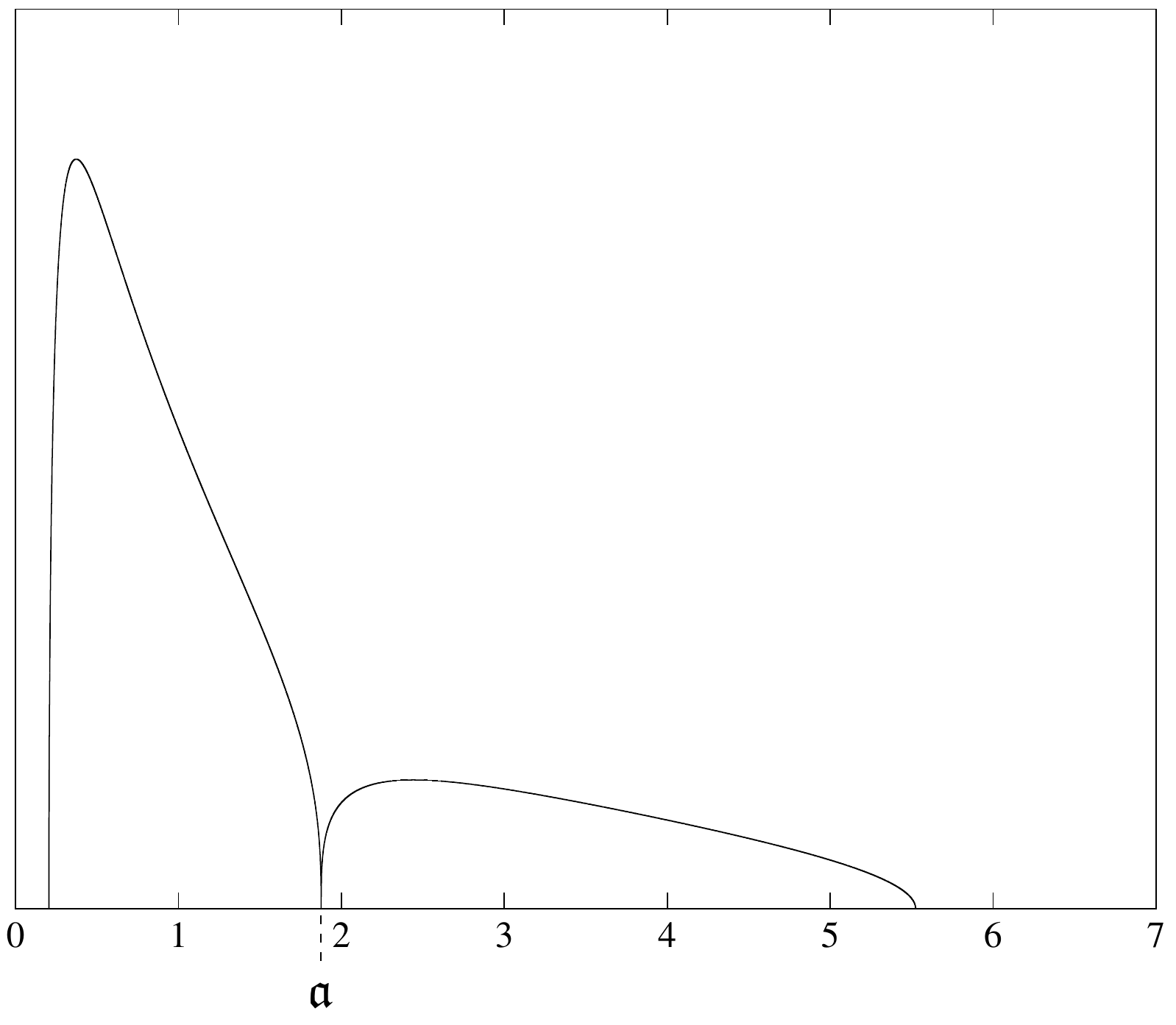}
\caption{Plot of the density of $\mu$ in the framework of 
Figure~\ref{fig:gcusp}.} 
\label{fig:cusp}
\end{figure}

\begin{proposition}
\label{prop global C} Let $\frak a$ be a cusp point,  set $\frak c = m(\frak a)$ and assume $\frak c\in D$. Then
\eq
\label{g cusp}
g(\frak c)=\frak a,\qquad g'(\frak c)=0,\qquad g''(\frak c)=0,
\qquad \text{and} \; \ g'''(\frak c)>0. 
\qe
Moreover, 
\eq
\label{sqrt cusp}
\rho(x)=\frac{\sqrt 3}{2\pi}\left(\,\frac{6}{g'''(\frak c)}\,\right)^{1/3}\big|x-\frak a\big|^{1/3}\;(1+o(1))\ , \qquad x\to\frak a\ .
\qe
Conversely, if $\frak c\in D$ satisfies $g'(\frak c)=g''(\frak c)=0$, then 
the real number $\frak a = g(\frak c) $ is a cusp point, $g'''(\frak c)>0$ 
and~\eqref{sqrt cusp}  holds true.
\end{proposition}
Propositions \ref{prop global 0} and \ref{prop global C} appear 
in~\cite{HHN2}.  Proposition~\ref{prop global C} is illustrated in Figures~\ref{fig:gcusp} 
and~\ref{fig:cusp}. 

\section{Local behavior}
\label{section local}

The study of the eigenvalues local behavior of random matrices is a central
topic in random matrix theory.  
When dealing with large Hermitian random matrices, it is recognized 
that the local correlation of the eigenvalues around an edge where
the density vanishes like a square root should be described 
by a particular point process involving the Airy kernel (see below), whose maximal particle's  distribution is known as the Tracy-Widom law. For instance, this has been established for unitary invariant random matrices and for Wigner matrices as well, see e.g. the surveys \cite{Dei,Erd} and references therein. Similarly, the
Bessel kernel is expected to describe the fluctuations around an hard edge
where the density vanishes like an inverse square root, and the Pearcey kernel
around a cusp with cubic root behavior.  Let us also mention that the sine
kernel is expected around a point where the density is positive, and that more
sophisticated behaviors have also been observed in matrix models where the
density vanishes like a rational power of different order, but we will not
further investigate these aspects here.  Another interesting feature not covered by this survey is the study of the random eigenvectors, see  e.g. \cite{BGN, BKYY}.

The purpose of this section is to present the central results of
\cite{HHN,HHN2}, where such typical local behaviors arise for the complex
correlated Wishart matrices under consideration at every edges and cusp points
satisfying a certain regularity condition. 

In fact, the free parameter family $(\nu_N)$ may have a deep impact on the limiting local fluctuations and one may not recover the expected fluctuations without further conditions. A first manifestation of this phenomenon is the Baik-Ben Arous-P\'ech\'e (BBP) phase transition, that we present in Section \ref{BBP}.
In a nutshell, this phase transition yields that slight variations on the family $(\nu_N)$ may modify the fluctuations at a soft edge and  may no longer be described by the Tracy-Widom law. Such phenomenas motivate the introduction of a regularity condition which essentially rules out this kind of behaviors.

In Section \ref{csq RC}, we provide the existence of finite $N$ approximations of the edges or cusp points under study and satisfying the regularity condition; the reader not interested in these precise definitions may skip this section. 

Next, in Section \ref{section Airy} we introduce the Airy kernel, the Tracy-Widom law, and state our results concerning the soft edges. In Section \ref{section Bessel} , we introduce the Bessel kernel and describe the fluctuations at the hard edge.  As an application, we provide in Section \ref{section condition} a precise description for the asymptotic  behavior of the condition number of $\bv M_N$. Finally, in Section \ref{section Pearcey} we introduce the Pearcey kernel, and state our result concerning the asymptotic behavior near a cusp point.

\subsection{The BBP phase transition and the regularity assumption}
\label{BBP}

First, assume $\bv \Sigma_N$ is the identity matrix, so that the (limiting) spectral distribution  $\nu$ of $\bv \Sigma_N$ is $\delta_1$, and hence the limiting density $\rho(x)$  is provided by~\eqref{dens-mp}.  
If $x_{\max}$ stands for the maximal eigenvalue of $\bv M_N$, then it has been established that $x_{\max}$ converges a.s. towards the right edge and fluctuates at the scale $N^{2/3}$ according to the Tracy-Widom law \cite{Jo}. Next, following Baik, Ben Arous and P\'ech\'e \cite{BBP},  assume instead $\bv \Sigma_N$ is a finite rank additive perturbation of the identity, meaning that the rank of the perturbation is independent on $N$. Thus we still have $\nu=\delta_1$ and the limiting density $\rho(x)$ remains unchanged. They established that if the strength of the perturbation is limited, then the behavior for $x_{\max}$ is the same as in the non-perturbed case, see \cite[Theorem 1.1(a), $k=0$]{BBP}. On the contrary, if the perturbation is strong enough, then $x_{\max}$ converges a.s.~outside of the bulk and the fluctuations are of different nature, see \cite[Theorem 1.1(b)]{BBP}. But in this case, one can consider instead the largest eigenvalue that actually converges to the right edge and show that the Tracy-Widom fluctuations still occur (this is a consequence of Theorem~\ref{th:fluctuations-TW} below). However, they also established there is an intermediary regime, where $x_{\max}$ converges a.s. to the right edge and fluctuates at the scale $N^{2/3}$ but not according to the Tracy-Widom law, see \cite[Theorem 1.1(a), $k>1$]{BBP}, hence leaving the random matrix universality class since the right edge exhibits a square root behavior.  Here the fluctuations are actually described by a deformation of the Tracy-Widom law, but in the general $\bv\Sigma_N$'s setting much exotic behaviors must be expected. 

In conclusion, although the eigenvalues global behavior only depends on the limiting parameters $\nu$ and $\gamma$, the local behavior is quite sensitive in addition to the mode of convergence of the spectral measure $\nu_N$ of $\bv\Sigma_N$ to its limit $\nu$. In order to obtain universal fluctuations in the more general setting under investigation, it is thus necessary to  add an extra condition for the edges, and actually for the cusp points too.  A more precise consideration of the non-universal intermediary regime considered by Baik, Ben Arous and P\'ech\'e reveals that, if we write the right edge as $g(\frak c)$, see Proposition \ref{prop global R} and the comments below, then some of the inverse eigenvalues of $\bv\Sigma_N$ converge towards $\frak c$. Recalling the $\lambda_j$'s stand for the eigenvalues of $\bv\Sigma_N$, this motivates us to introduce the following condition.

\begin{definition} A real number  $\frak c$ satisfies the \textbf{regularity condition} if 
\eq
\label{RC}
\liminf_{N\to\infty}\min_{j=1}^n\left|\frak c-\frac{1}{\lambda_j}\right|>0.
\qe
Moreover, if $\frak c$ satisfies the regularity condition, we then say that $g(\frak c)$ is \textbf{regular}.
\end{definition} 

\begin{remark}
\label{reg edge} 
 Propositions \ref{prop global R} and \ref{prop global L} tell us that every soft edge reads $g(\frak c)$ for some $\frak c\in \overline D$. In fact, since $g(0)=+\infty$ and $\mathrm{Supp}(\mu)$ is compact, necessarily $\frak c\neq 0$. If we moreover assume the soft edge to be regular then, since  by definition $D = \{ x \in \R : x \neq 0,\, x^{-1} \not\in \supp(\nu) \}$ and because $\nu_N$ converges weakly to $\nu$, necessarily $\frak c\in D$. In particular, Propositions \ref{prop global R} and \ref{prop global L} yield that at a regular soft edge the density show a  square root behavior.  As regards the hard edge, the analogue of the regularity condition turns out to be $\liminf_N\lambda_1>0$ and is therefore contained in Assumption~\ref{ass:nu}.


\end{remark}

\begin{remark}
We show in \cite{HHN} that, if $\gamma>1$, then the leftmost edge $\frak a$ is always regular. Namely there exists $\frak c\in D$ which is regular such that $\frak a= g(\frak c)$. In fact, we have $\frak c<0$.
\end{remark}

Before we state our results on the eigenvalues local behavior around the regular edges or the cusp points, we now provide the existence of the appropriate scaling parameters we shall use in the later statements. 

\subsection{Consequences of the regularity condition and finite $N$ approximations for the edges and the cusp points}
\label{csq RC}

Recall from Section~\ref{section global} that the Cauchy transform of the 
limiting eigenvalue distribution $\mu$ of $\bv M_N$ is defined as the unique
solution $m \in \C_-$ of the fixed-point equation~\eqref{cauchy eq}. 
We now consider the probability measure $\mu_N$ induced after replacing  $(\gamma, \nu)$ by its finite horizon analogue 
$(n/N, \nu_N)$ in this equation (we recall $\nu_N$ was introduced in \eqref{nuN}). Namely, let 
$\mu_N$ be the probability measure whose Cauchy transform is defined as 
the unique solution $m \in \C_-$ of the fixed-point equation 
\[
m = \left( z - \frac nN \int \frac{\lambda}{1 - m \lambda} \nu_N(\d\lambda) 
\right)^{-1} \ . 
\]
The probability measure $\mu_N$ should be thought of as a deterministic 
approximation of the distribution of the eigenvalues of $\bv M_N$ at finite 
$N$, and is referred to as the \textbf{deterministic equivalent} of the 
spectral measure of $\bv M_N$. The measure $\mu_N$ reads 
\[
\mu_N(\d x)=\left(1-\frac{n}{N}\right)^+ \delta_0 +  \rho_N(x)\d x,
\]
and one can apply all the results stated in Section \ref{section global} to 
describe $\rho_N$, after replacement of $g$ with 
\eq
\label{gN}
g_N(z)=\frac 1{z}+\frac nN\int\frac{\lambda}{1-z\lambda}\nu_N(\d \lambda) .
\qe

Recalling $D$ has been introduced in \eqref{D}, the following proposition encodes the essential consequence of the regularity condition.

\begin{proposition} 
\label{gN->g}If $\frak c\in D$ satisfies the regularity condition \eqref{RC}, then there exists $\delta>0$  such that $g_N$ is holomorphic on $\{z\in\C:\;|z-\frak c|< \delta\}\subset D$ for every $N$ large enough and converges uniformly towards $g$ there.
\end{proposition}

It is an easy consequence of Montel's theorem. Now, if a sequence of holomorphic function $h_N$ converges uniformly to a (holomorphic) function $h$ on an open disc, then a standard result from complex analysis provides that the $k$th order derivative $h_N^{(k)}$  also converges uniformly to $h^{(k)}$ there, for every $k\geq 1$. Moreover, Hurwitz's theorem  states that, if $h$ has a zero $\frak c$ of multiplicity $\ell$ in this disc, then $h_N$ has exactly $\ell$ zeros, including multiplicity,  converging towards $\frak c$ as $N\to\infty$. 

Thus, as a consequence of the previous proposition, by applying Hurwitz's 
theorem to $g_N'$ (and the symmetry $g_N'(\bar z)=g_N'( z)$), it is  easy to obtain the following statement. 

\begin{proposition} 
\label{edge cN}
Assume $\frak c\in D$ satisfies the regularity condition \eqref{RC} and moreover 
\[
g'(\frak c)=0,\qquad g''(\frak c)<0,\qquad resp.\quad  g''(\frak c)>0.
\]
Then there exists a sequence $(\frak c_N)$, unique up to a finite number of terms, converging to $\frak c$ and such that, for every $N$ large enough, we have $\frak c_N\in D$ and
\[
\lim_{N\to\infty}g_N(\frak c_N)=g(\frak c),\qquad g_N'(\frak c_N)=0,\qquad g_N''(\frak c_N)<0,\qquad resp.\quad  g_N''(\frak c_N)>0.
\]
\end{proposition}

 Having in mind  Propositions \ref{prop global L} and \ref{prop global R},   this proposition thus states that if one considers a regular left (resp. right) soft edge $\frak a$, and thus $\frak a=g(\frak c)$ with $\frak c\in D$ by Remark \ref{reg edge}, then there exists a sequence, unique up to a finite number of terms, of  left (resp. right) soft edges $\frak a_N=g_N(\frak c_N)$ for the deterministic equivalent $\mu_N$ converging towards $\frak a$. These soft edges $(\frak a_N)$ are finite $N$ approximations of the edge $\frak a$, while the $\frak c_N$'s are finite $N$ approximations of the preimage $\frak c$.

When dealing with regular cusp points, the situation is slightly more delicate. The reason for this is that if $\frak a=g(\frak c)$ is a regular cusp point, then $\frak c$ is now a zero of multiplicity two for $g'$. By applying Hurwitz's theorem to $g'_N$ as above, one would obtain two sequences of non-necessarily real zeros for $g_N'$ converging towards $\frak c$. It is actually more convenient to apply Hurwitz's theorem to $g_N''$ instead, in order to get the following statement. 

\begin{proposition} 
\label{cusp cN}
Assume $\frak c\in D$ satisfies the regularity condition \eqref{RC} and moreover 
\[
g'(\frak c)=0,\qquad g''(\frak c)=0,\qquad \text{(hence} \ g'''(\frak c)>0 
\ \text{by Prop. }\ref{prop global C} \text{).}  
\]
Then there exists a sequence $(\frak c_N)$, unique up to a finite number of terms, converging to $\frak c$ and such that, for every $N$ large enough, we have $\frak c_N\in D$ and
\[
\lim_{N\to\infty}g_N(\frak c_N)=g(\frak c),\qquad \lim_{N\to\infty}g_N'(\frak c_N)=0,\qquad g_N''(\frak c_N)=0,\qquad  g_N'''(\frak c_N)>0.
\]
\end{proposition}

Notice that Proposition \ref{cusp cN} doesn't guarantee that 
$g_N'(\frak c_N)=0$. Hence, a cusp point is not necessarily the limit of cusp points of the deterministic equivalents $\mu_N$. As we shall see in Section~\ref{section Pearcey}, the speed at which $g_N'(\frak c_N)$ goes to zero will actually influence the local behavior around the cusp.

\begin{definition} \label{def cN}
Given a soft left edge, resp. right edge, resp. cusp point  $\frak a$  which is regular, and thus $\frak a=g(\frak c)$ with $g'(\frak c)=0$ and $g''(\frak c)<0$, resp. $g''(\frak c)>0$, resp. $g''(\frak c)=0$ and  $g'''(\frak c)>0$,  the \textbf{sequence associated with $\frak a$} is the sequence $(\frak c_N)$ provided by Propositions \ref{edge cN} and \ref{cusp cN}.    

\end{definition}

Equipped with Propositions \ref{edge cN} and \ref{cusp cN}, we are now in position to state the results concerning the local asymptotics.

\subsection{The Airy kernel and Tracy-Widom fluctuations at a soft edge}
\label{section Airy}

Given a  function $\K(x,y)$ from $\R\times\R$ to $\R$ satisfying appropriate conditions, one can consider its associated determinantal point process,  which is a simple point process  on $\R$ having as correlation functions the determinants $\det[\,\K(y_i,y_j)]$. More precisely, it is a probability distribution $\p$ over the configurations $(y_i)$ of real numbers (the particles), namely over discrete subsets of $\R$ which are locally finite, characterized in the following way: For every $k\geq 1$ and any test function $\Phi:\R^k\to\R$, 
\[
\mathbb E\left[   \sum_{y_{i_1}\neq \, \cdots \,\neq \, y_{i_k}}  \Phi(y_{i_1},\ldots,y_{i_k})\right]=\int_\R \cdots\int_{\R} \Phi(y_1,\ldots,y_k)\det\Big[\K(y_i,y_j)\Big]_{i,j=1}^k\d y_1\cdots \d y_k,
\] 
where the sum runs over the $k$-tuples of pairwise distinct particles of the configuration $(y_i)$.  Hence the correlation between the particles $y_i$'s is completely encoded by the kernel $\K(x,y)$. In particular,  the inclusion-exclusion principle yields a closed formula for the gap probabilities in terms of Fredholm determinants. Namely, for any interval $J\subset \R$, the probability that no particle lies in $J$ reads
\[
\p\Big( (y_i)\cap J =\emptyset\Big)= 1+\sum_{k=1}^\infty\frac{(-1)^k}{k!}\int_J \cdots\int_{J}\det\Big[\K(y_i,y_j)\Big]_{i,j=1}^k\d y_1\cdots \d y_k\, ,
\]
and the latter is the Fredholm determinant $\det(I-\K)_{L^2(J)}$ of the integral operator acting on $L^2(J)$ with kernel $\K(x,y)$, provided it makes sense.  We refer to \cite{Hu,JoR} for further information on determinantal point processes. 

Consider the \textbf{Airy point process} $\p_{\Ai}$ which is defined as the determinantal point process on $\R$ associated with the Airy kernel
\eq
\label{Kai}
\K_{\Ai}(x,y)=\frac{{\rm Ai}(x){\rm Ai}'(y)-{\rm Ai}(y){\rm Ai}'(x)}{x-y},
\qe
where the Airy function 
\[
\mathrm {Ai}(x)=\frac1\pi\int_0^{\infty} \cos\left(\frac{u^3}{3}+ux\right)\d u
\]
is a solution of  the differential equation $ f''(x)=xf(x)$.

The configurations $(y_i)$ generated by the Airy point process  a.s. involve an infinite number of particles but have a largest particle $y_{\max}$. The distribution of  $y_{\max}$ is the \textbf{Tracy-Widom law} (see e.g. \cite[Section 2.2]{JoR}), and its distribution function reads, for every $s\in\R$,
\eq
\label{TW def}
\p_\Ai\big (y_{\max}\leq  s \big )=\p_\Ai\Big ((y_{i})\cap (s,+\infty)=\emptyset \Big )=\det(I-\K_{\Ai})_{L^2(s,\infty)}.
\qe
Tracy and Widom \cite{TW1} established the  famous representation
\[
\p_\Ai\big (y_{\max}\leq  s \big )=\exp\left({-\int_{s}^\infty}(x-s)q(x)^2\d x\right),
\]
where $q$ is the Hastings-McLeod solution of the Painlev\'e II equation, namely the unique solution of  $f''(x)=f(x)^3+xf(x)$ with boundary condition $f(x)\sim {\rm Ai}(x)$ as $ x\rightarrow+\infty$.\\

Recalling that $g_N$ has been introduced in \eqref{gN}, we are now in position
to describe the eigenvalues local behavior around regular soft edges. In the
three upcoming theorems, we denote by 
$\tilde x_1\leq \cdots \leq \tilde x_{n}$ the ordered eigenvalues of 
$\widetilde{\bv M}_N$. We also use the notational convention $\tilde x_0 = 0$ 
and $\tilde x_{n+1} = + \infty$.

\begin{theorem} \label{th:fluctuations-TW}  
Let $\frak a$ be a right edge and assume it is regular. Writing 
$\frak a = g(\frak c)$, let 
$\phi(N) = \max\{ j : \lambda_j^{-1} > \frak c\}$. Then, almost surely, 
\eq
\label{extr R}
 \tilde x_{\phi(N)}\xrightarrow [N\to\infty]{} \frak a,\qquad 
\liminf_{N\to\infty} \big( \tilde x_{\phi(N)+1} - \frak a\big)>0.
\qe
Moreover, let  $(\frak c_N)_N$ be the sequence associated with $\frak a$ as in  Definition \ref{def cN}. Set 
\eq
\label{cst R}
\frak a_N=g_N(\frak c_N),\qquad \sigma_N=\left(\frac{2}{g_N''(\frak c_N)}\right)^{1/3},
\qe
so that $\frak a_N\to \frak a$, $\frak c_N\to \frak c$, and 
$\sigma_N\to (2/g''(\frak c))^{1/3}>0$ as $N\to\infty$. Then, for every 
$s\in\R$,
\eq
\label{TW right}
\lim_{N\rightarrow\infty} \p\Big( 
  N^{2/3}\sigma_N\big( \tilde x_{\phi(N)}-\frak a_N\big)\leq s\Big)
=\p_\Ai\big (y_{\max}\leq  s \big ).
\qe
\end{theorem}

\begin{remark} Let us stress that the sequence $\phi(N)$ may be non-trivial even when considering the rightmost edge: As explained in Section \ref{BBP},  it is indeed possible that a certain amount  of eigenvalues (possibly infinite) will converge outside of the limiting support. Thus, if we assume the rightmost edge is regular, then Theorem \ref{th:fluctuations-TW}   states that there exists a maximal eigenvalue $x_{\phi(N)}$ which actually converge to the rightmost edge and fluctuates according to the Tracy-Widom law.

\end{remark}
Let us comment  on the history of this theorem. The Tracy-Widom fluctuations
have been first obtained by Johansson \cite{Jo} for the maximal eigenvalue when
$\bv\Sigma_N$ is the identity.   Baik, Ben Arous and P\'ech\'e \cite{BBP} then
proved this  still holds true when $\bv\Sigma_N$ is a finite rank perturbation
of the identity, provided the perturbation is small enough.  Assuming a
condition which is equivalent to the regularity condition \eqref{RC} and that
the maximal eigenvalue converges towards the rightmost edge, El Karoui
\cite{EK} established the Tracy-Widom fluctuations for the maximal eigenvalue
for general $\bv\Sigma_N$'s assuming $\gamma\leq 1$, and Onatski \cite{On} got
rid of the  last restriction. The statement on the existence of extremal
eigenvalues converging to each regular right edge is \cite[Theorem 2]{HHN}  and
essentially relies on the exact separation results of Bai and Silverstein
\cite{BS1,BS2}; the definition of the sequence $\phi(N)$ indexing these
extremal eigenvalues relies on these results. Finally, the Tracy-Widom
fluctuations for the extremal eigenvalues associated to any right regular edge
is \cite[Theorem 3-(b)]{HHN}.

We now provide a similar statement for the left soft edges.

\begin{theorem} 
\label{th:fluctuations-TW2} 
Let $\frak a$ be a left edge of the bulk. 
If $\gamma > 1$ and $\frak a$ is the leftmost edge of the bulk, set 
$\phi(N) = n-N+1$. Otherwise, assume that $\frak a>0$ is regular, write 
$\frak a = g(\frak c)$ and set 
$\phi(N) = \min \{ \lambda_j : \lambda_j^{-1} < \frak c \}$. Then, almost
surely, 
\eq
\label{extr L}
 \tilde x_{\phi(N)}\xrightarrow [N\to\infty]{} \frak a,\qquad 
\liminf_{N\to\infty} \big( \frak a-\tilde x_{\phi(N)-1} \big)>0.
\qe
Moreover, let  $(\frak c_N)_N$ be the sequence associated with $\frak a$ as in
Definition~\ref{def cN}. Set 
\eq
\label{cst L}
\frak a_N=g_N(\frak c_N),\qquad \sigma_N=\left(\frac{2}{-g_N''(\frak c_N)}\right)^{1/3},
\qe
so that $\frak a_N\to \frak a$, $\frak c_N\to \frak c$, and 
$\sigma_N\to (-2/g''(\frak c))^{1/3}>0$ as $N\to\infty$. Then, for every 
$s\in\R$,
\eq
\label{TW left}
\lim_{N\rightarrow\infty} \p\Big( 
  N^{2/3}\sigma_N\big( \frak a_N-\tilde x_{\phi(N)}\big)\leq s\Big)
=\p_\Ai\big (y_{\max}\leq  s \big ).
\qe
\end{theorem}

Prior to this result, which is a combination of Theorem 2 and Theorem 3-(a) from \cite{HHN}, the Tracy-Widom fluctuations for the smallest random eigenvalue when $\bv\Sigma_N$ is the identity has been obtained by Borodin and Forrester \cite{BF}.   

Let us also mention that when $\nu$ is the sum of two Dirac masses, a local uniform convergence to the Airy kernel (which is a weaker statement than the Tracy-Widom fluctuations) at every (regular) right and soft left edges follows from \cite{LW}, see also \cite{Mo1,Mo2}.

Finally, we state our last result, concerning the asymptotic independence of the Tracy-Widom fluctuations at a finite number of regular soft edges. For a more precise statement, we refer to  \cite[Theorem 4]{HHN}.
 
\begin{theorem} 
\label{th:independence} 
Let $(\frak a_j)_{j\in J}$ be a finite collection of soft edges, and assume all
these edges are regular. For each $j\in J$, consider the rescaled eigenvalue
$N^{2/3}\sigma_{N,j}(\tilde x_{\phi_j(N)}-\frak a_{N,j})$ associated with the soft
edge $\frak a_j$ provided by \eqref{extr R}--\eqref{cst R} and \eqref{extr
L}--\eqref{cst L}.  Then the family of random variables
$\{N^{2/3}\sigma_{N,j}(\tilde x_{\phi_j(N)}-\frak a_{N,j})\}_{j\in J}$ becomes
asymptotically independent as $N\to\infty$.
\end{theorem}

The asymptotic independence has  been previously established for the smallest and largest  eigenvalues when $\bv\Sigma_N$ is the identity by Basor, Chen and Zhang \cite{BCZ}.

\begin{remark}
\label{KY univ} The results presented in this survey rely on the fact that the entries of $\bv X_N$ are  complex Gaussian random variables, 
a key assumption in order to take advantage of the determinantal structure of the eigenvalues of the model under study. A recent work \cite{knowles-yin-2014-preprint} by 
Knowles and Yin enables to transfer the results of Theorems \ref{th:fluctuations-TW}, \ref{th:fluctuations-TW2} and \ref{th:independence}
to the case of  complex, but not necessarily Gaussian, random variables. Indeed,  by combining the local convergence to the limiting distribution established in \cite{knowles-yin-2014-preprint} together with Theorems \ref{th:fluctuations-TW}, \ref{th:fluctuations-TW2} and 
\ref{th:independence}, one obtains Tracy-Widom fluctuations and asymptotic independence in this more general setting, provided that the entries of matrix $\bv{X}_N$ fulfill some moment condition. Let us stress that the case of 
real Gaussian random variables (except the largest one covered in \cite{lee-schnelli-preprint}), of important interest in statistical applications, remains open.
\end{remark}

We now turn to the hard edge and the Bessel point process.

\subsection{The Bessel point process at the hard edge}
\label{section Bessel}

The \textbf{Bessel point process} $\p_\Be^{(\alpha)}$ of parameter $\alpha\in\mathbb Z$ is the determinantal point process on $\R_+$ associated with the kernel 
\eq
\label{Bessel kernel}
\K_{\Be}^{(\alpha)}(x,y)=\frac{\sqrt{y}\, J_\alpha(\sqrt{x})J_\alpha'(\sqrt{y})-\sqrt{x}\, J_\alpha'(\sqrt{x})J_\alpha(\sqrt{y})}{2(x-y)} ,
\qe
where the Bessel function of the first kind $J_\alpha$ with parameter $\alpha$ is defined for $x\geq0$ by 
\eq
\label{series rep Bessel}
J_\alpha(x) = \left( \frac x2\right)^\alpha \sum_{n=0}^\infty \frac{(-1)^n}{n!\,  \Gamma(n+\alpha+1)} \left( \frac x2\right)^{2n}
\qe
and satisfies the differential equation $x^2 f''(x)+xf'(x)+(x^2-\alpha^2)f(x)=0$.  

The configurations $(y_i)$ generated by the Bessel point process  a.s. have an infinite number of particles $y_i$  but have a smallest particle $y_{\min}$. The law of  $y_{\min}$ is characterized, for every $s>0$, by
\eq
\label{hard TW def}
\p_\Be^{(\alpha)}\big (y_{\min}\geq   s \big )=\p_\Be^{(\alpha)}\Big ((y_{i})\cap (0,s)=\emptyset \Big )=\det(I-\K_{\Be}^{(\alpha)})_{L^2(0,s)}.
\qe
When $\alpha=0$, this reduces to an exponential law of parameter $1$, namely  $\p_\Be^{(0)}\big (y_{\min}\geq   s \big )=e^{-s}$, as observed by Edelman \cite{Ed}. In the general case, Tracy and Widom obtained the representation  \cite{TW2},
\[
\p_\Be^{(\alpha)}\big (y_{\min}\geq   s \big )=\exp\left(-\frac{1}{4}\int_0^s(\log s -\log x) q(x)^2\d x\right),
\]
where $q$ is the solution of a differential equation which is reducible to a particular case of the Painlev\'e V equation (involving $\alpha$ in its parameters) and boundary condition $q(x)\sim J_\alpha(\sqrt x)$ as $x\to 0_+$. 

Recalling the $\lambda_j$'s are the eigenvalues of $\bv\Sigma_N$ and their distributional limit $\nu$ has a compact support in $(0,\infty)$, we now provide our statement concerning the eigenvalues local behavior around the hard edge.

\begin{theorem}
\label{th Bessel}
Assume that $n=N+\alpha$ with $\alpha\in\mathbb Z$ independent of $N$ and set
\[
\sigma_N=-2g_N''(\infty)=\frac{4}{N}\sum_{j={1}}^n\frac{1}{\lambda_j} \ ,\qquad \zeta_N=-\frac{4}{3}g'''_N(\infty)=\frac{8}{N}\sum_{j={1}}^n\frac{1}{\lambda_j^2}\ .
\] 
Thus $\sigma_N\to 4\int\lambda^{-1}\nu(\d \lambda)>0$ and $\zeta_N\to 8\int\lambda^{-2}\nu(\d \lambda)>0$ as $N\to\infty$. 

Let $x_{\min}$ be the smallest random eigenvalue of $\bv M_N$. Then, for every $s>0$, we have 
\eq
\label{Bessel first}
\lim_{N\to\infty}\p\Big(N^2\sigma_N\,x_{\min}\geq s\Big)=\p_\Be^{(\alpha)}\big (y_{\min}\geq   s \big ).
\qe
Furthermore, we have the expansion as $N\to\infty$,
\eq
\label{Bessel next}
\p\Big(N^2\sigma_N\,x_{\min}\geq s\Big)=\p_\Be^{(\alpha)}\big (y_{\min}\geq   s \big )-\frac1N\left(\frac{\alpha \zeta_N}{\sigma_N^2}\right) \, s\frac{\d }{\d s}\p_\Be^{(\alpha)}\big (y_{\min}\geq   s \big )+O\left(\frac{1}{N^2}\right)\, .
\qe
\end{theorem}

When $\bv\Sigma_N$ is the identity, the convergence \eqref{Bessel first} has been established by Forrester \cite{F}. As for the next order term \eqref{Bessel next}, it has been obtained when $\bv\Sigma_N$ is the identity by Perret and Schehr \cite{PS} and Bornemann \cite{bornemann-2014-note}, motivated by a question raised by Edelman, Guionnet and P\'ech\'e in \cite{EGP}. The statement \eqref{Bessel first} has been first obtained in \cite{HHN}, while the stronger statement \eqref{Bessel next} is \cite[Theorem 6]{HHN2}.

\subsection{Application to condition numbers}
\label{section condition}

In this subsection, we study the fluctuations of the ratio 
$$
\kappa_N=\frac {x_{\max}}{x_{\min}} 
$$
of the largest to the smallest random eigenvalue of $\bv M_N$. Notice that 
if $n \geq N$, then $\kappa_N$ is the condition number of $\bv M_N$ while if 
$n \leq N$, then $\kappa_N$ is the condition number of $\widetilde{\bv M}_N$. 
The condition number is a central object of study in numerical linear 
algebra~\cite{von-neumann-goldstine-47,von-neumann-goldstine-51}. 
Using our previous results, we can obtain an asymptotic description for 
$\kappa_N$.  
Let us emphasize that the leftmost edge $\frak a$ of the support of $\rho$ is positive 
if and only if $\gamma \neq 1$, see \cite[Proposition 3]{HHN2}. 

\begin{proposition} \label{prop:condition-number}
Assume $\gamma \neq 1$. Denote by $\frak a$ the leftmost edge and by $\frak b$
the rightmost one. Assume that $\frak a,\frak b$ are regular, $x_{\min}\to\frak a$ and $x_{\max}\to \frak b$ 
a.s. as $N\to\infty$ (that $\frak a $ is regular and $x_{\min}\to\frak a$ a.s. is always true when $\gamma >1$). Write $\frak a=g(\frak c)$, $\frak b=g(\frak d)$, consider the sequences $(\frak c_N)$ and $(\frak d_N)$ associated with $\frak c$ and $\frak d$ respectively (see Definition \ref{def cN}) and set
\[
\frak a_N=g_N(\frak c_N)\ ,\quad \sigma=\left(\frac{2}{-g''(\frak c)}\right)^{1/3},\quad \qquad 
\frak b_N=g_N(\frak d_N)\ ,\quad \delta=\left(\frac{2}{g''(\frak d)}\right)^{1/3}\ .
\]
Then, 
$$
\kappa_N \xrightarrow[N\to \infty]{a.s.} \frac{\frak b}{\frak a}\qquad \textrm{and}\qquad 
N^{2/3} \left( 
\kappa_N - \frac{\frak b_N}{\frak a_N}
\right) \xrightarrow[N\to \infty]{\mathcal D} \frac {X} {\delta\frak a} + \frac {\frak b Y}{\sigma \frak a^2}
$$
where $\xrightarrow[]{\mathcal D}$ stands for the convergence in distribution 
and where $X$ and $Y$ are two independent random variables with the 
Tracy-Widom distribution.
\end{proposition}

We now handle the case where $\gamma=1$.

\begin{proposition}\label{prop:condition-number-gamma-1}
Assume $n=N+\alpha$, where $\alpha\in \mathbb{Z}$ is independent of $N$, and moreover $x_{\max}\to \frak b$ a.s. for some $\frak b>0$.  Then,  
$$
\frac 1{N^2} \kappa_N \xrightarrow[N\to \infty]{\mathcal D} \frac {4\frak b } {X} \left(\int \lambda^{-1}\nu(\d \lambda)\right)
$$
where $\p(X\geq s)=\p_\Be^{(\alpha)}\big (y_{\min}\geq   s \big )$ for every $s>0$.
\end{proposition}

\begin{remark} 
Interestingly, in the square case where $\gamma=1$, the fluctuations of the
largest eigenvalue $x_{\max}$ have no influence on the 
fluctuations of $\kappa_N$ as these are imposed by the limiting distribution 
of $x_{\min}$ and the a.s. limit $\frak b$ of $x_{\max}$.
\end{remark}

Finally, we turn to the eigenvalues local behavior near a cusp point.

\subsection{The Pearcey kernel at a cusp point}
\label{section Pearcey}
Given any $\tau\in\R$, following \cite{TW3} we  introduce the Pearcey-like integral functions
\[
\phi(x)=\frac{1}{2i\pi}\oint_{\Sigma} e^{xz-\tau z^2/2+z^4/4} \d z,\qquad \psi(y)=\frac{1}{2i\pi}\int_{-i\infty}^{i\infty} e^{-yw+\tau w^2/2-w^4/4} \d w,
\]
where the contour $\Sigma$ consists in two rays going from $\pm e^{i\pi/4}\infty$ to zero, and two rays going from zero to $\pm e^{-i\pi/4}$. They satisfy the respective differential equations 
\[
\phi'''(x)-\tau\phi'(x)+x\phi(x)=0, \qquad \psi'''(y)-\tau\psi'(y)-y\psi(y)=0.
\]
The \textbf{Pearcey point process} $\p_\Pe^{(\tau)}$ is the determinantal point process associated with the Pearcey kernel
\eq
\label{KPe}
\K_\Pe^{(\tau)}(x,y)=\frac{\phi''(x)\psi(y)-\phi'(x)\psi'(y)+\phi(x)\psi''(y) -\tau\psi(x)\psi(y)}{x-y}.
\qe
This process has been first introduced by Br\'ezin and Hikami \cite{BH1,BH2} when $\tau=0$, and subsequent generalizations have been considered by Tracy and Widom \cite{TW3}.

The configurations $(y_i)$ generated by the Pearcey point process are a.s. infinite and do not have a largest nor smallest particle. With this respect, the  quantities  of interest here are the gap probabilities of the Pearcey point process, defined for every $0<s<t$ by
\eq
\label{def gap Pearcey}
\p_{\Pe}^{(\tau)}\Big( (y_i)\cap [s,t]=\emptyset\Big)=\det(I-\K_{\Pe}^{(\tau)})_{L^2(s,t)} .
\qe
Seen as a function of $s$, $t$ and $\tau$, the $\log$ of the righthand side of \eqref{def gap Pearcey} is know to satisfy a system of PDEs, see \cite{TW3, BC12, ACvM12}.

In \cite[Theorem 5]{HHN2}, we prove the following statement.

\begin{theorem} \label{main cor} 
Let $\frak a = g(\frak c)$ be a cusp point such that $\frak c \in D$, and 
assume it is regular. Let $(\frak c_N)$ be the sequence associated with 
$\frak a$ as in Definition~\ref{def cN}. Assume moreover
the following decay assumption holds true: There exists $\kappa\in\R$ such that
\eq
\label{speed asump}
\sqrt N\, g_N'(\frak c_N)\xrightarrow[N\to\infty]{} \kappa \ .
\qe
We set
\eq
\label{constants}
\frak a_N=g_N(\frak c_N),\qquad \sigma_N=\left(\frac{6}{g_N^{(3)}(\frak c_N)}\right)^{1/4},\qquad \tau=-\kappa\left(\frac{6}{ g^{(3)}(\frak c)}\right)^{1/2},
\qe
so that $\frak a_N\to\frak a$ and $\sigma_N\to\left(6/g'''(\frak c)\right)^{1/4}>0$ as $N\to\infty$. Then, for every $s>0$, we have 
\eq
\label{gap conv}
\lim_{N\to\infty}\p\Big( \big(N^{3/4}\sigma_N(x_i-\frak a_N)\big)\cap[-s,s]=\emptyset\Big) = \p_{\Pe}^{(\tau)}\Big( (y_i)\cap [-s,s]=\emptyset\Big),
\qe
where the $x_i$'s are the random eigenvalues of $\bv M_N$.
\end{theorem}

\noindent This result has been obtained by Mo when $\bs\Sigma_N$ has exactly two distinct eigenvalues  \cite{Mo2}.

As advocated in Section \ref{csq RC}, the precise decay for $g_N'(\frak c_N)\to 0$ does influence the eigenvalues local behavior near a cusp  (see Proposition \ref{cusp cN} and the discussion below).  Our assumption \eqref{speed asump} covers  the general case where $\sqrt N g_N'(\frak c_N)\to0$, and hence the limiting kernel is $\K_\Pe^{(0)}(x,y)$ introduced by Br\'ezin and Hikami, and the limiting regime where $\sqrt N g_N'(\frak c_N)$ has a limit as well. 

\begin{figure}[h]
\centering
\includegraphics[width=0.7\linewidth]{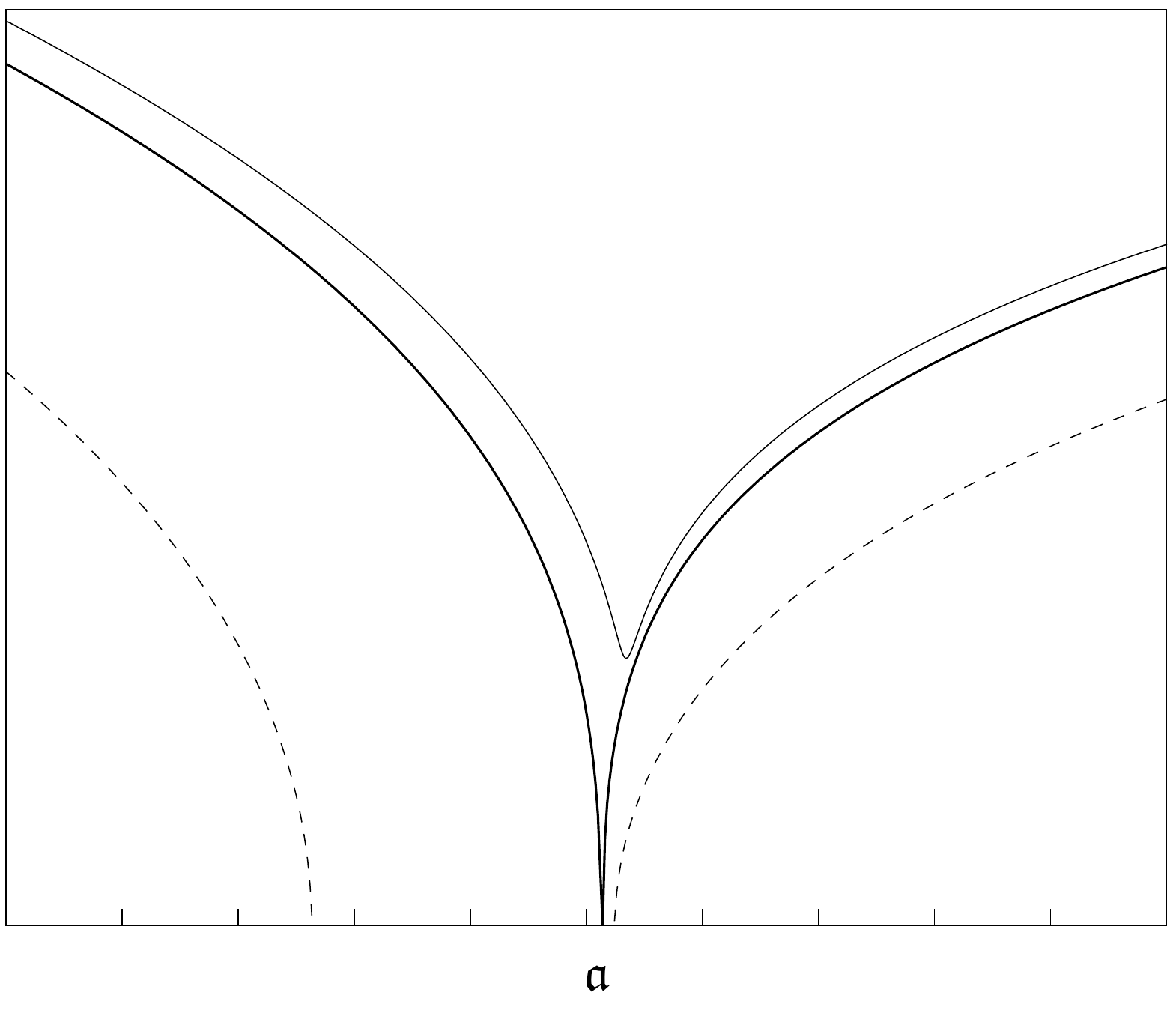}
\caption{Zoom of the density of $\mu_N$ near the cusp point
$\frak a$. The thick curve is the density of $\mu$ in the framework
of Figure~\ref{fig:cusp}. The thin curve (resp.~the dotted curve) is the 
density of $\mu_N$ when $\sqrt{N} g'_N(\frak c_N) > 0$
(resp.~$\sqrt{N} g'_N(\frak c_N) < 0$).}  
\label{fig:zoomcusp}
\end{figure}

\begin{remark}{\bf (erosion of a valley)}
In the case where this limit $\kappa$ in \eqref{speed asump} is positive, the deterministic equivalent measure $\mu_N$ will not feature a cusp but rather a valley that will become deeper
as $N\to\infty$, see the thin curve in Figure \ref{fig:zoomcusp}. The density of $\mu_N$ will always be positive near the cusp and the condition 
$$
g'_N(\frak c_N) \sim \frac{\kappa}{\sqrt{N}}
$$  
should be thought of as a speed condition of the erosion of the valley. 
\end{remark}

\begin{remark}{\bf (moving cliffs)} 
In the case where $\kappa<0$ in \eqref{speed asump}, $g'_N(\frak c_N)$ is
always negative for $N$ large enough. In particular, there exists a small
$N$-neighborhood of $\frak c_N$ whose image by $g_N$ is outside the support of
$\mu_N$: There is a small hole in the support of $\mu_N$ but the two connected
components move towards one another (moving cliffs), see the dotted curve in
Figure \ref{fig:zoomcusp}. In this case, the condition 
$$ 
g'_N(\frak c_N) \sim \frac{\kappa}{\sqrt{N}} 
$$ can also be interpreted as a speed condition at which the
cliffs approach one another.  
\end{remark}

\begin{remark}{\bf (slow decay)}
The slow decay setting where $$\sqrt N g_N'(\frak c_N)\to \pm\infty$$ is not covered by our results. In this case, we do not expect the Pearcey point process to arise anymore, and refer to Section \ref{OP} for further discussion.
\end{remark}

\section{Sketches of the proofs}\label{sec:proofs}

In this section,  we provide an outline for the proofs of the results presented in Section \ref{section local}.

\subsection{The random eigenvalues of $\bv M_N$ form a determinantal point process }

The key input on which all our proofs are based on, is that when the elements
of $\bv X_N$ are complex Gaussian (Assumption~\ref{ass:gauss}), the configuration
of the random eigenvalues $x_i$'s of $\bv M_N$ form a determinantal point
process with an explicit kernel. More precisely Baik, Ben Arous and P\'ech\'e
provided in \cite{BBP} a formula for that kernel, to which they give credit to
Johansson. It is given by the following double complex integral 
\eq
\label{KN}
\K_N(x,y)=\frac{ N}{(2i\pi)^2}\oint_{\Gamma}\d z\oint_{\Theta} \d w\,\frac{1}{w-z} e^{- Nx(z-\frak q) +Ny(w-\frak q)}\left(\frac{z}{w}\right)^{ N}\prod_{j=1}^n\left(\frac{w-\lambda_j^{-1}}{z-\lambda_j^{-1}}\right),
\qe
 where the $\frak q\in\R$ is a free parameter (see \cite[Remark 4.3]{HHN}) and we recall the $\lambda_j$'s are the eigenvalues of $\bv\Sigma_N$.    
 The contours $\Gamma$ and $\Theta$ are disjoint and closed, both oriented counterclockwise, such that $\Gamma$ encloses all the $\lambda_j^{-1}$'s  whereas $\Theta$ encloses the origin. 
\begin{figure}[ht]
\centering 
\includegraphics[width=0.7\linewidth]{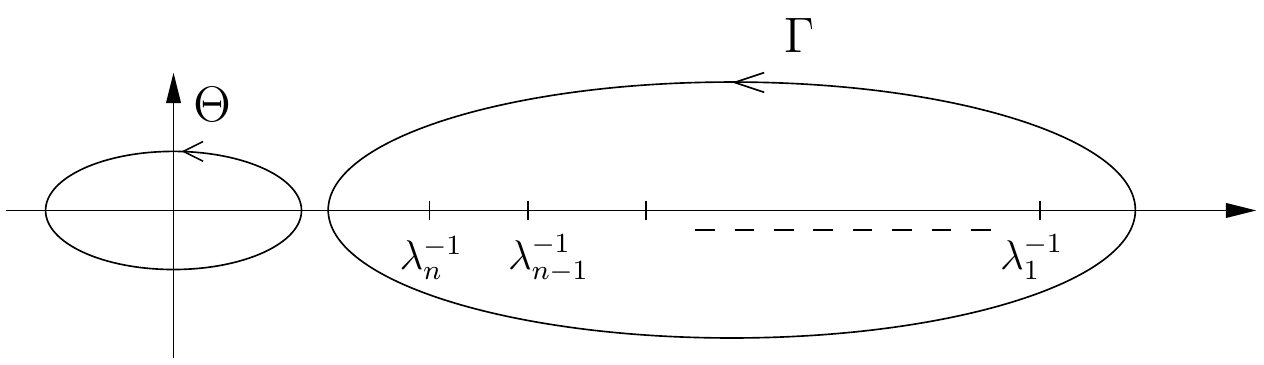}
\caption{The contours of integration} 
\label{fig:paths} 
\end{figure} 
  
  \begin{remark} The main ingredient to obtain this determinantal representation is the Harish-Chandra-Itzykson-Zuber integral formula, which allows to write a particular integral over the unitary group in terms of determinants, see \cite[Section 2.1]{BBP}. The analogue of this integral formula does not seem to exist for correlated Wishart matrices with real or quaternionic entries, and thus the determinantal structure seems only available in the complex setting.   
  
  \end{remark}
  If we consider a random configuration of the form
  \[
  \big(N^\beta\sigma_N(x_i-\frak a_N)\big)\ , 
 \]
 then a change of variables yields it is a determinantal point process with kernel
 \eq
 \label{res KN}
 \frac{1}{N^\beta\sigma_N}\K_N\left( \frak a_N+\frac{x}{N^\beta\sigma_N}, \frak a_N+\frac{y}{N^\beta\sigma_N}\right) ,
 \qe
 where $\K_N$ is as in \eqref{KN}. Hence, the study of the eigenvalues local behavior boils down to the asymptotic analysis as $N\to\infty$ of kernels of the form \eqref{res KN} with different choices for the scaling parameters $\beta,\sigma_N,\frak a_N$.

 \subsection{Modes of convergence}
 
In order to prove the convergence \eqref{Bessel first} at the hard edge,  it is enough to establish a local uniform convergence on $\R_+\times\R_+$  for the kernel \eqref{res KN} to the Bessel kernel $\K_{\Be}^{(\alpha)}(x,y)$, after choosing appropriately the scaling  parameters. Similarly, the local uniform convergence on $\R\times\R$ to the Pearcey kernel  $\K_{\Pe}^{(\tau)}(x,y)$ yields the convergence \eqref{gap conv} for the gap probabilities around a cusp. The convergences \eqref{TW right} and \eqref{TW left} to the Tracy-Widom law however require a stronger mode of convergence  (such as the trace-class norm convergence, or the Hilbert-Schmidt norm plus trace convergence, for the associated operators acting on $L^2(s,\infty)$, for every $s\in\R$; we refer to \cite[Section 4.2]{HHN} for further information).  This essentially amounts to obtain  a local uniform convergence on $(s,+\infty)\times (s,+\infty)$ plus tail estimates for $\K_N(x,y)$.

From now, we shall disregard these convergence issues and provide heuristics on why the Airy kernel, the Bessel kernel  and the Pearcey kernel 
should appear in different scaling limits.

\subsection{Towards the Airy kernel}

Here we provide an heuristic for the convergence to the Airy kernel. The gap to be filled in order to make this sketch of  a proof mathematically rigorous can be found in \cite{HHN};  this  heuristic may actually serve as a roadmap for the quite lengthy and technical proof we provided there.

\begin{figure}[ht]
\centering 
\includegraphics[width=0.4\linewidth]{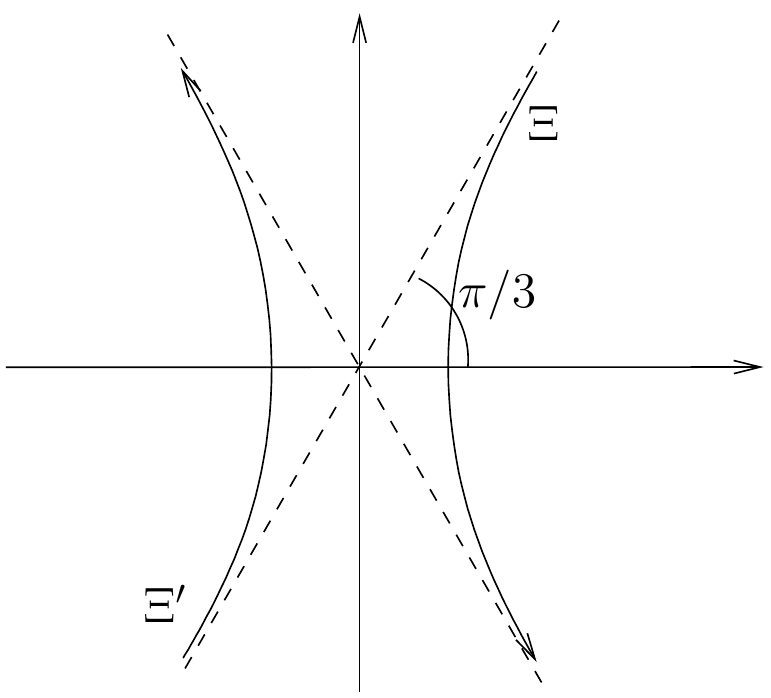}
\caption{The paths of integration for the Airy kernel} 
\label{fig:paths-Airy} 
\end{figure}

Since we are dealing with contours integrals, it is more convenient to use the following alternative representation for the Airy kernel \eqref{Kai},
\eq
\label{Airy cont}
\K_\Ai(x,y)=
\frac{1}{(2i\pi)^2}\int_{\Xi}\d z\int_{\Xi'} \d w \,\frac{1}{w-z} e^{-xz+yw+z^3/3-w^3/3},
\qe
which is based on the contour integral formula for the Airy function (see e.g.
the proof of \cite[Lemma 4.15]{HHN}). The contours $\Xi$ and $\Xi'$ are
disjoint and unbounded contours, and $\Xi$ goes from $e^{i\pi/3}\infty$ to
$e^{-i\pi/3}\infty$ whereas $\Xi'$ goes from $e^{-2i\pi/3}\infty$ to
$e^{2i\pi/3}\infty$, as shown on Figure~\ref{fig:paths-Airy}. 

Consider the scaling parameters associated with a soft edge provided in Theorem \ref{th:fluctuations-TW}; we thus focus on the right edge setting, but the situation for a left edge  is similar. More precisely, by using the formula \eqref{KN} where we take $\frak q=\frak c_N$, we investigate 
\begin{align}
\label{KN airy int}
 & \frac{1}{N^{2/3}\sigma_N}\K_N\left(\frak a_N+ \frac{x}{N^{2/3}\sigma_N},\frak a_N+\frac{y}{N^{2/3}\sigma_N}\right)\\
  & \quad = \quad\frac{N^{1/3}}{(2i\pi)^2\sigma_N}\oint_\Gamma\d z\oint_\Theta \d w\, \frac{1}{w-z} e^{-N^{1/3}x(z-\frak c_N)/\sigma_N+N^{1/3}y(w-\frak c_N)/\sigma_N} e^{Nf_N(z)-Nf_N(w)} \nonumber,
\end{align}
where we introduced the map
\eq
\label{fN}
f_N(z)=-\frak a_{N}(z-\frak  c_N)+\log(z)-\frac{1}{N}\sum_{j=1}^n\log(1-\lambda_j z).
\qe
After performing  the change of variables $z\mapsto \frak c_N+\sigma_N z/N^{1/3}$ and $w\mapsto \frak c_N+\sigma_N w/N^{1/3}$,  the right-hand side of \eqref{KN airy int} becomes
\begin{multline}
\label{int change airy}
\frac{1}{(2i\pi)^2}\oint_{\phi_N(\Gamma)}\d z\oint_{\phi_N(\Theta)} \d w\, \frac{1}{w-z} e^{-xz+yw+N f_N\big(\frak c_N + \sigma_N\frac{z}{N^{1/3}}\big)-N f_N\big(\frak c_N + \sigma_N\frac{w}{N^{1/3}}\big)},
\end{multline}
where we set for convenience $\phi_N(z)= N^{1/3}(z-\frak c_N)/\sigma_N$.

Next, recalling $g_N$ was introduced in \eqref{gN}, the crucial observation that
\eq
\label{fgN0}
f_N'(z)= g_N(z)-\frak a_N
\qe 
allows to infer  on the local behavior of $f_N$ around $\frak c_N$. More precisely, since by definition of the scaling parameters we have
\[
g_N(\frak c_N)=\frak a_N,\qquad g_N'(\frak c_N)=0,\qquad g_N''(\frak c_N)\xrightarrow[N\to\infty]{}g''(\frak c)>0,
\]
a Taylor expansion for $f_N$ around $\frak c_N$ yields the approximation
\eq
\label{preapprox airy}
N\Big (f_N(\frak c_N + \sigma_N\frac{z}{N^{1/3}})-f_N(\frak c_N)\Big) \; \simeq  \; \frac{1}{3}z^3\ ,
\qe
 where the constant $1/3$ comes from the definition of $\sigma_N$, see \eqref{cst R}. In conclusion, after plugging \eqref{preapprox airy} into \eqref{int change airy}, we obtain the approximation 
\begin{multline}
\label{KN airy int 2}
  \frac{1}{N^{2/3}\sigma_N}\K_N\left(\frak a_N+ \frac{x}{N^{2/3}\sigma_N},\frak a_N+\frac{y}{N^{2/3}\sigma_N}\right)\\
    \simeq \frac{1}{(2i\pi)^2}\oint_{\phi_N(\Gamma)}\d z\oint_{\phi_N(\Theta)} \d w\, \frac{1}{w-z} e^{-xz+yw+ z^3/3-w^3/3},
\end{multline}
and we can almost read the Airy kernel \eqref{Airy cont}, up to contour deformations.

To frame the previous heuristic into a rigorous mathematical setting, a few technical points should be addressed, since of course the approximation \eqref{preapprox airy} is only valid when $|z|$ is not too large and the contours appearing in the Airy kernel are unbounded. With this respect, the standard move is to split  the contours $\Gamma$ and $\Theta$ into different parts and then to deform each part in an appropriate way.

In a neighborhood of $\frak c_N$, after simple transformations, one chooses $\Gamma$ and $\Theta$ to match with the contours of the Airy kernel there, and then justify rigorously the approximation \eqref{KN airy int 2} after restriction of $z,w$ to that neighborhood.  This can be done by quantifying the approximation \eqref{preapprox airy} and then performing  tedious but rather simple computations. 

Then, outside of this neighborhood, one proves that the remaining of the integrals don't contribute in the large $N$ limit. In the present setting of a general matrix $\bv \Sigma_N$, this is the hard part of the proof. To do so, one  establishes the existence of admissible deformations for the contours $\Gamma$ and $\Theta$ so that they complete the Airy contours truncated on a neighborhood of $\frak c_N$, and where the contribution coming from the term $\exp\{Nf_N(z)-Nf_N(w)\}$ brings an exponential decay on the remaining part. This can be done by looking for the so-called steepest descent/ascent contours (i.e. contours on which $\re f_N$ is decreasing/increasing), and this was the strategy used by Baik, Ben Arous, P\'ech\'e  \cite{BBP} and El Karoui  \cite{EK} when dealing with the rightmost edge. When considering any right or left soft edge,  following this strategy  requires to consider many sub-cases  and to perform  again most of the computations in several case. In \cite{HHN}, we instead developed  a unified (abstract) method to provide the existence of appropriate contours by mean of the maximum principle for subharmonic functions.

For the reader interested in having a  look at the proofs of \cite{HHN}, let us mention  it turns out it is more convenient to work at a scale where the contours $\Gamma,\Theta$ live in a bounded domain, and this is the reason why we did not performed there the changes of variables $z\mapsto \frak c_N+\sigma_N z/N^{1/3}$ and $w\mapsto \frak c_N+\sigma_N w/N^{1/3}$ as we did in the present heuristic.\\

\subsection{Towards the Pearcey kernel}

Now we turn to the heuristics for Pearcey kernel and refer to \cite{HHN2} for a rigorous proof.  The setting is essentially the same as in the Airy case, except that now $\frak c_N$ is a simple zero for $g_N''$ instead of $g_N'$, and $g_N'(\frak c_N)\to 0$,  which entails a different behavior  for the map $f_N$ near $\frak c_N$.

We start with the alternative representation   for the Pearcey kernel \eqref{KPe},
\eq
\label{Pearcey cont}
\K_{\Pe}^{(\tau)}(x,y)=
\frac{1}{(2i\pi)^2}\int_{\Xi}\d z\int_{-i\infty}^{\, i\infty} \d w \,\frac{1}{w-z} e^{-xz-\frac{\tau z^2}2+\frac{z^4}4+yw+\frac{\tau w^2}2-\frac{w^4}4},
\qe
 where the contour $\Xi$ is disjoint to the imaginary axis, and has two components. The first part goes from $e^{i\pi/4}\infty$ to $e^{-i\pi/4}\infty$, whereas the other part goes from $e^{-3i\pi/4}\infty$ to $e^{3i\pi/4}\infty$. See \cite{TW3} for a proof (and also \cite{BH1} when $\tau=0$). Notice also the symmetry $\K_{\Pe}^{(\tau)}(x,y)=\K_{\Pe}^{(\tau)}(-x,-y)$ which follows from the change of variables $z,w\mapsto -z,-w$.

Consider the scaling parameters associated with a regular cusp point provided in Theorem~\ref{main cor}. By using the formula \eqref{KN} where we choose $\frak q=\frak c_N$, we now consider 
\begin{align}
\label{KN Pearcey int}
 & \frac{1}{N^{3/4}\sigma_N}\K_N\left(\frak a_N+ \frac{x}{N^{3/4}\sigma_N},\frak a_N+\frac{y}{N^{3/4}\sigma_N}\right)\\
  & \quad = \quad\frac{N^{1/4}}{(2i\pi)^2\sigma_N}\oint_\Gamma\d z\oint_\Theta \d w\, \frac{1}{w-z} e^{-N^{1/4}x\frac{(z-\frak c_N)}{\sigma_N}+N^{1/4}y\frac{(w-\frak c_N)}{\sigma_N}} e^{Nf_N(z)-Nf_N(w)} \nonumber,
\end{align} where the map $f_N$ is the same as in \eqref{fN}. After  the change of variables $z\mapsto \frak c_N+\sigma_N z/N^{1/4}$ and $w\mapsto \frak c_N+\sigma_N w/N^{1/4}$,  the right-hand side of \eqref{KN Pearcey int} reads
\begin{multline}
\label{int change Pearcey}
\frac{1}{(2i\pi)^2}\oint_{\phi_N(\Gamma)}\d z\oint_{\phi_N(\Theta)} \d w\, \frac{1}{w-z} e^{-xz+yw+N f_N\big(\frak c_N + \sigma_N\frac{z}{N^{1/4}}\big)-N f_N\big(\frak c_N + \sigma_N\frac{w}{N^{1/4}}\big)},
\end{multline}
where we introduced for convenience $\phi_N(z)= N^{1/4}(z-\frak c_N)/\sigma_N$.

In this setting, the definition of the scaling parameters yields
\[
g_N(\frak c_N)=\frak a_N,\qquad \sqrt N g'_N(\frak c_N)\xrightarrow[N\to\infty]{}\kappa, \qquad g_N''(\frak c_N)=0,\qquad g_N'''(\frak c_N)\xrightarrow[N\to\infty]{}g'''(\frak c)>0.
\]
Recalling the identity \eqref{fgN0} and the definition \eqref{constants} of $\tau$, a Taylor expansion around $\frak c_N$ then yields the approximation
\eq
\label{preapprox Pearcey}
N\Big (f_N(\frak c_N + \sigma_N\frac{z}{N^{1/4}})-f_N(\frak c_N)\Big) \; \simeq  \; -\frac{\tau}{2}z^2+\frac{1}{4}z^4\ ,
\qe
where the constant $1/4$ comes from the definition of $\sigma_N$, see \eqref{constants}. Thus, by plugging \eqref{preapprox Pearcey} into \eqref{int change Pearcey}, we obtain the approximation 
\begin{multline}
\label{KN Pearcey int 2}
  \frac{1}{N^{3/4}\sigma_N}\K_N\left(\frak a_N+ \frac{x}{N^{3/4}\sigma_N},\frak a_N+\frac{y}{N^{3/4}\sigma_N}\right)\\
    \simeq \frac{1}{(2i\pi)^2}\oint_{\phi_N(\Gamma)}\d z\oint_{\phi_N(\Theta)} \d w\, \frac{1}{w-z} e^{-xz-\frac{\tau z^2}2+\frac{z^4}4+yw+\frac{\tau w^2}2- \frac{w^4}4},
\end{multline}
and we can almost see the Pearcey kernel \eqref{Pearcey cont}, up to contour deformations.

As to make this approximation rigorous, the method is the same as for the Airy kernel. Let us mention that the abstract argument we mention previously for the existence of appropriate contour deformations also applies  in this setting.

\subsection{Towards the Bessel kernel}

Finally, we provide heuristics for the appearance of the Bessel kernel and refer to \cite{HHN2} for a rigorous proof.  The main input here is, according to Section \ref{section global}, the critical point $\frak c$ associated with the hard edge is now located at  infinity.

The first step is to write the Bessel kernel \eqref{Bessel kernel} as the double contour integral,
\eq
\label{Bessel cont}
\K_{\Be}^{(\alpha)}(x,y)=
\frac{1}{(2i\pi)^2}\left(\frac{y}{x}\right)^{\alpha/2}\oint_{|z|=\, r} \frac{\d z}{z}\oint _{|w|=\, R}\frac{\d w}{w} \,\frac{1}{z-w}\left(\frac{z}{w}\right)^\alpha e^{-\frac xz+\frac z4+\frac yw-\frac w4},
\qe
where $0<r<R$, and which is provided in \cite[Lemma 6.2]{HHN}. The contours of integration are circles oriented counterclockwise. Let us stress this formula is only available when $\alpha\in\mathbb Z$, since otherwise the term $(z/w)^\alpha$ in the integrand would not make sense on the whole of the integration contours.

Setting $\sigma_N$ as in Theorem \ref{th Bessel} and using the formula \eqref{KN} where we choose $\frak q=0$, we now consider 
\begin{align}
\label{KN Bessel int}
 & \frac{1}{N^{2}\sigma_N}\K_N\left(\frac{x}{N^{2}\sigma_N},\frac{y}{N^{2}\sigma_N}\right)\\
  & \quad = \; \frac{1}{(2i\pi)^2N\sigma_N}\oint_{\Gamma} \d z \oint_{\Theta} \d w\,\frac{1}{w-z} \left(\frac{z}{w}\right)^{N}e^{- \frac{zx}{N\sigma_N}+\frac{wy}{N\sigma_N}}\prod_{j=1}^n\frac{w-\lambda_j^{-1}}{z-\lambda_j^{-1}} \nonumber.
\end{align} 
Having in mind the critical point is located at infinity, we  perform the change of variables $z\mapsto N\sigma_N/z$ and $w\mapsto N\sigma_N/w$, so that  the right-hand side of \eqref{KN Bessel int} reads
\begin{multline}
\label{int change Bessel}
\frac{1}{(2i\pi)^2}\oint_{\phi_N(\Gamma)}\frac{\d z}z
\oint_{\phi_N(\Theta)} \frac{\d w}w\, \frac{1}{z-w} 
\left(\frac{z}{w}\right)^\alpha 
e^{-\frac xz+\frac yw-N G_N(z)+ NG_N(w)},
\end{multline}
where we introduced the maps
\[
G_N(z)=\frac{1}{N}\sum_{j=1}^n\log\left(\frac{z}{N\sigma_N} -\lambda_j\right)
\]
and  $\phi_N(z)= N\sigma_N/z$. We emphasize that, during the previous step, we used that $n=N+\alpha$ and witnessed a cancellation leading to the term $(z/w)^\alpha$, which does not depend on $N$. 

Now, a Taylor expansion of $G_N$ around zero yields the approximation
\eq
\label{preapprox Bessel}
N\Big (G_N(z)-G_N(0)\Big) \; \simeq  \; - \frac{z}{4N}\ .
\qe
Thus, by plugging \eqref{preapprox Bessel} into \eqref{int change Bessel}, we obtain the approximation 
\begin{multline}
\label{KN Bessel int 2}
   \frac{1}{N^{2}\sigma_N}\K_N\left(\frac{x}{N^{2}\sigma_N},\frac{y}{N^{2}\sigma_N}\right)\\
    \simeq \frac{1}{(2i\pi)^2}\oint_{\phi_N(\Gamma)}\frac{\d z}z\oint_{\phi_N(\Theta)} \frac{\d w}w\, \frac{1}{z-w} \left(\frac{z}{w}\right)^\alpha e^{-\frac xz+\frac yw+\frac z4-\frac w4},
\end{multline}
and we can almost see the Bessel kernel \eqref{Bessel cont}, up to contour deformations and the prefactor $(y/x)^{\alpha/2}$. Finally, in order to deal with that prefactor,  one considers the operator $\E$ of multiplication by $x^{\alpha/2}$ acting on $L^2(0,s)$, and then use that replacing the Bessel kernel $\K_\Be^{(\alpha)}(x,y)$ by the kernel of the operator $\E\K_\Be^{(\alpha)}\E^{-1}$, which is \eqref{Bessel cont} without the prefactor $(y/x)^{\alpha/2}$, leaves the Fredholm determinant $\det(I-\K_\Be^{(\alpha)})_{L^2(0,s)}$ invariant. 

To make this heuristic rigorous, the method is far less demanding than in the
setting of the Airy or the Pearcy kernel. Indeed, in the present setting one
can  legally deform the contours $\Gamma$ and $\Theta$ in such a way that
$\phi_N(\Gamma)$ and $\phi_N(\Theta)$ match the integration contours for the
Bessel kernel \eqref{Bessel cont}. After making this deformation, 
a simple Taylor expansion of the map $G_N$ around zero will be enough to 
establish the convergence towards the Bessel kernel and therefore to 
obtain~\eqref{Bessel first}. 

By pushing the Taylor expansion  \eqref{preapprox Bessel} one step further, 
one can also obtain the more accurate estimate~\eqref{Bessel next}, provided with an identity  involving the resolvent of the Bessel kernel  established by Tracy and Widom \cite{TW2}. We refer the reader to \cite{HHN2} for further information. 

\section{Open questions}
\label{OP}

The results presented here naturally entail a number of open questions that we list below.

\begin{enumerate}
%
\item {\bf At the edge of the definition domain and exotic local behaviors}. The results on the  eigenvalues local behavior presented in this survey only concern edges or cusp points $\frak a$ which read $\frak a=g(\frak c)$ with $\frak c\in D$. If we focus on the rightmost edge for the sake of simplicity, then Proposition \ref{prop global R} states  this edge may actually belong to $g(\partial D)$ (notice that this cannot happen if the limiting spectral distribution $\nu$ of $\bv\Sigma_N$ is a finite combination of Dirac measures). In this case, the square root behavior of the density around this edge is not guaranteed  anymore, and the laws describing the fluctuations of the eigenvalues near such an edge seem completely unknown and a priori different from Tracy-Widom distribution.  We believe  the fluctuations will actually depend on $\nu$, and hence lie outside of the random matrix universality class. Quite interestingly, the same phenomenon arise in the study of the additive deformation of a GUE random matrix \cite{capitaine-peche-2014-preprint} and random Gelfand-Tsetlin patterns \cite{DM1,DM2}.

\item {\bf Alternative regime at a cusp point I}. In the context of Theorem \ref{main cor}, our speed assumption \eqref{speed asump} does not cover the following case
$$
\sqrt N\, g_N'(\frak c_N)\xrightarrow[N\to\infty]{} +\infty \ .
$$
This condition corresponds to the situation where the density of  the deterministic equivalent $\mu_N$ is positive in a neighborhood of $\frak c_N$. 
It essentially states that the bulk of $\mu_N$ will degenerate into a cusp around $g(\frak c)$ quite slowly and we do not expect to witness  
Pearcey-like fluctuations around $g_N(\frak c_N)$ anymore. We believe instead that the sine kernel will arise at the scale $\sqrt{N/g_N'(\frak c_N)}$, which strictly lies  in between $N^{1/2}$ and $N^{3/4}$. 
\item {\bf Alternative regime at a cusp point II}. Another case that is not covered by our assumption \eqref{speed asump} is when
$$
\sqrt N\, g_N'(\frak c_N)\xrightarrow[N\to\infty]{} -\infty \ .
$$
In this case, $\frak c_N$ lies outside the support of  $\mu_N$ and $g'_N$ has two distinct real zeroes near $\frak c_N$, say $\frak c_{N,1}$ and $\frak c_{N,2}$.
Hence, for $N$ sufficiently large, $g_N(\frak c_{N,1})$ and  $g_N(\frak c_{N,2})$ both correspond to  edges of the support of $\mu_N$ which both converge towards the cusp point  $g(\frak c)$. The previous condition entails that such a convergence will happen at a quite slow rate. In this case, we do not expect to observe the Pearcey kernel around $g_N(\frak c_N)$ either, because of the absence of particles, but we believe a local analysis around the edge $g_N(\frak c_{N,1})$ or $g_N(\frak c_{N,2})$ may uncover the Airy kernel at an intermediary scale. 

\item {\bf Study of the fluctuations at the hard edge in more general cases}.
The hard edge fluctuations were described  here when 
$n = N + \alpha$ with a fixed $\alpha\in \mathbb Z$ , but the hard edge is always present as soon as $n/N\to1$. Thus it would be of interest to describe the hard edge fluctuations in  more
general situations, for example when $\alpha=\alpha(n) \to +\infty$ so that $n/N \to 1$.  In this case one would expect Tracy-Widom fluctuations near the leftmost edge $\frak a_N$ of $\mu_N$, the latter being positive but converging to zero as $N\to \infty$.

\item {\bf Non-Gaussian entries}. 

All the fluctuations results presented here rely on the fact that the entries
of matrix $\bv X_N$ are complex Gaussian. It is however of major interest, for
applications and for the general theory as well, to study the universality of
such results for non-Gaussian complex random variables. As explained in Remark
\ref{KY univ},  the Tracy-Widom  fluctuations for the extremal eigenvalues
associated with any regular soft edges are now established in this general
setting (under some moment conditions for the entries), by combining theorems
\ref{th:fluctuations-TW} and \ref{th:fluctuations-TW2} with Knowles and Yin's
recent preprint \cite{knowles-yin-2014-preprint} (see also \cite{bao-et-al}).
However,  natural related questions remain open: Would it be
possible to describe the fluctuations at

\begin{itemize}
\item[(a)]  the hard edge, for general complex entries?
\item[(b)]  a (regular) cusp point, for general complex entries?
\end{itemize}

Another universality class of interest is the case where the matrix $\bv X_N$ has real entries. In this case, the techniques based on the determinantal structure of the 
eigenvalues are no longer available. Lee and Schnelli \cite{lee-schnelli-preprint} recently succeeded to establish GOE Tracy-Widom  fluctuations of the largest eigenvalue when the entries  are real Gaussian or simply real (with subexponential decay), under the assumption that the covariance matrix $\bv \Sigma_N$ is diagonal. The techniques developed by  Knowles and Yin \cite{knowles-yin-2014-preprint} enable to relax the diagonal assumption for the covariance matrix 
$\bv \Sigma_N$. A number of questions remain open: Would it be possible to describe the fluctuations at
\begin{itemize}
\item[(c)] any (regular) soft edge, when the entries are real Gaussian?
\item[(d)] the hard edge, when the entries are real (Gaussian or not)?
\item[(e)] a (regular) cusp point, when the entries are real (Gaussian or not)?
\end{itemize}

\end{enumerate}

\end{document}